\documentclass[11pt]{amsart}
\usepackage[utf8]{inputenc}
\usepackage[english]{babel}

\usepackage{amscd,amssymb}
\usepackage{amsmath}
\usepackage{amsthm}
\usepackage{graphicx}
\usepackage{fancybox}
\usepackage{epic,eepic}
\usepackage{amstext}
\usepackage{mathtools}
\usepackage{esint}
\usepackage[square,numbers]{natbib}
\usepackage{booktabs}
\usepackage{hyperref}
\usepackage[noabbrev, capitalise, nameinlink]{cleveref}
\crefname{equation}{}{}
\crefname{assumption}{Assumption}{Assumptions}
\crefformat{equation}{\textup{#2(#1)#3}}
\usepackage{comment}
\usepackage{mathtools}
\usepackage{tikz,pgfplots}
\usepackage{subcaption}

\def\V{V}

\def\W{\mathcal{W}}

\def\Pnull{\mathbb{P}^0}

\def\elem{T}

\def\T{\mathcal{T}}
\def\TH{\mathcal{T}_H}

\def\C{\mathcal{C}}

\def\PiH {\Pi_H}

\def\one{\mathbf{1}}

\newcommand{\mnormf}[1]{{| #1 |}_{M}}
\newcommand{\mnormfo}[2]{{| #1 |}_{M,#2}}

\newcommand{\lnormf}[1]{{| #1 |}_{L}}

\newcommand{\vsp}[2]{{\left( #1,\,#2 \right)}}
\newcommand{\vspf}[2]{{( #1,\,#2 )}}

\newcommand{\with}{\,:\,}

\definecolor{myBlue}{RGB}{113,104,238} 
\definecolor{myGreen}{RGB}{114,175,30} 
\definecolor{myRed}{RGB}{180,50,50}  
\definecolor{myOrange}{RGB}{225,92,22}

\newtheorem{theorem}{Theorem}[section]

\newtheorem{lemma}[theorem]{Lemma}
\newtheorem{assumption}[theorem]{Assumption}

\theoremstyle{definition}

\theoremstyle{remark}
\newtheorem{remark}[theorem]{Remark}
\numberwithin{theorem}{section}
\numberwithin{equation}{section}
\numberwithin{table}{section}
\numberwithin{figure}{section}

\allowdisplaybreaks

\begin{document}
	
\title[An algebraic multiscale method]{An algebraic multiscale method for spatial network models}
\author[M.~Hauck, R.~Maier, A.~M\aa lqvist]{Moritz Hauck$^\dagger$, Roland Maier$^{\ddagger}$, Axel M\aa lqvist$^\dagger$}
\address{${}^{\dagger}$ Department of Mathematical Sciences, Chalmers University of Technology and University of Gothenburg, 41296 Göteborg, Sweden}
\email{hauck@chalmers.se,\, axel@chalmers.se}
\address{${}^{\ddagger}$ Institute for Applied and Numerical Mathematics, Karlsruhe Institute of Technology, Englerstr.~2, 76131 Karlsruhe, Germany}
\email{roland.maier@kit.edu}
\maketitle

\begin{abstract}
	In this work, we present a multiscale approach for the reliable coarse-scale approximation of spatial network models represented by a linear system of equations with respect to the nodes of a graph. The method is based on the ideas of the Localized Orthogonal Decomposition (LOD) strategy and is constructed in a fully algebraic way. This allows to apply the method to geometrically challenging objects such as corrugated cardboard. In particular, the method can also be applied to finite difference or finite element discretizations of elliptic partial differential equations, yielding an approximation with similar properties as the LOD in the continuous setting. We present a rigorous error analysis of the proposed method under suitable  assumptions on the network. Moreover, numerical examples are presented that underline our theoretical results.
\end{abstract}

\vspace{1cm}
\\[2ex]

\section{Introduction}

Spatial network models are a powerful tool for simplifying the representation of rather complex systems, such as blood vessels~\cite{FKOWW22} or pore formations in porous materials~\cite{CEPT12}, into one-dimensional structures. Another interesting application is the modeling of paper, where wooden fibers form a net-like structure~\cite{KMM20}.  
In this work, we aim to achieve reliable coarse-scale approximations of spatial networks. We consider a spatial network $\mathcal{G}=(\mathcal{N},\mathcal{E})$ consisting of nodes and edges together with an operator~$K$ defined on $\mathcal G$ which is spectrally equivalent to the graph Laplacian. 
The linear system which seeks the solution $u$ to
\begin{equation}
 	\label{eq:modelproblem0}
 	Ku=f
\end{equation}
with prescribed boundary conditions and a source term $f$, is typically poorly conditioned due to, e.g., complex geometries of the underlying network or highly heterogeneous material coefficients.  
Note that standard iterative solvers such as algebraic multigrid methods (see, e.g.,~\cite{LB12,XZ17,HWZ21}) may converge slowly when applied to \cref{eq:modelproblem0} because they often do not adequately handle the complex geometry of the spatial network.  
Alternatively, we may appropriately embed the network into a bounded domain, where we require that the network obeys a certain structure so that it essentially behaves like a continuous material at sufficiently coarse scales. In such a setting, one can apply  finite element techniques and  appropriate interpolation operators to solve~\cref{eq:modelproblem0} on an artificial coarse mesh. For example, multiscale methods can be used to obtain reliable coarse-scale approximations to the solution of~\cref{eq:modelproblem0}. A prominent example is the Localized Orthogonal Decomposition (LOD)  method \cite{MaP14,HeP13,MalP20,Peterseim2021,Mai20ppt}, which constructs problem-adapted coarse-scale approximation spaces with uniform approximation properties at the cost of a logarithmic overhead in the support of the basis functions when compared to classical finite elements. Recently, substantial progress has been made in \cite{HaPe21b} in further localizing the construction; see also~\cite{pumslod}. In conjunction with spatial networks, such approaches have been used in~\cite{KMM20,GoHeMa22,EGHKM24,HaM22}. Other multiscale approaches specifically designed for spatial networks are introduced in~\cite{EILRW09,ILW10}. 

In this paper, we generalize the above-mentioned LOD-based multiscale approaches for spatial networks to avoid the introduction of an auxiliary finite element mesh of the bounding domain. Instead, we use an algebraic approach where the discretization scales are defined using a suitable graph partitioning algorithm. We emulate the general construction of the LOD and in particular the bubble constructions from~\cite{HaPe21,DHM23}; see also the related technique known as gamblets~\cite{Owh17,OwhS19}. The resulting algebraic multiscale method provides reliable coarse approximations of the solution~$u$ to \eqref{eq:modelproblem0} and can be applied to complex geometries such as corrugated cardboard. While the network is still embedded into a physical domain, for example to define boundary conditions, our final algorithm is algebraic and avoids any direct relation between the network and the surrounding domain. We identify necessary assumptions on the network for the construction to work and present several examples to which the general methodology can be applied. In addition, a rigorous error analysis is presented and several numerical examples are given.

The remainder of this paper is organized as follows. In \cref{s:mp}, we introduce the concept of spatial networks as well as the considered model problem. Further, examples of problems which are included in the presented framework are given. Then, in \cref{sec:network}, abstract coarse partitions and necessary assumptions on the network are introduced. A prototypical LOD-based method, which achieves reliable coarse-scale approximations is introduced and analyzed in \cref{sec:protomethod}, and its decay and localization properties are discussed in \cref{sec:decloc}. Practical versions of the approximation are defined and analyzed in \cref{sec:prac} and studied numerically in \cref{sec:numexp}.

\section{A spatial network model}\label{s:mp}

We consider a connected graph $\mathcal{G} = (\mathcal{N}, \mathcal{E})$ of nodes
$\mathcal{N}$ and edges $\mathcal{E}$. For $x,\,y\in \mathcal{N}$, we write $x \sim y$ if $x$ and $y$ are connected by an edge. Therefore,
\begin{equation*}
	\mathcal{E} = \{ \{x, y\} \with x \sim y,\; x,y\in\mathcal N\}
\end{equation*}
and we assume $\mathcal{E}$ to be embedded into a domain $\Omega \subset \mathbb{R}^n$ for some $n \in \mathbb N$, which we assume to be scaled to unit size. Such graphs, which can be embedded in a spatial domain, are called \emph{spatial networks}.
We define the function space~$\hat V$ as the space of all real-valued functions supported on the set of nodes $\mathcal{N}$. Furthermore, we denote by
\begin{equation*}
V = \{ v \in \hat V\with v(x) = 0\text{ if } x \in \Gamma\}
\end{equation*} 
the subset of $\hat {V}$ satisfying homogeneous Dirichlet boundary conditions on the boundary segment $\Gamma  \subset \partial \Omega$. 
In the following, we will assume that $\mathcal N \cap \Gamma\neq \emptyset$, i.e., there exists at least one Dirichlet node. For any subset of nodes $S \subset \mathcal N$, we define the restrictions of $V$ and $\hat V$ to $S$ by $V_S$ and $\hat V_S$, respectively. 

\subsection{Model problem}

Before presenting the model problem of this work, we first introduce the mass and stiffness operators on the spatial network $\mathcal G$.  
Given $\alpha \in [0,2]$, we define for any $x \in \mathcal{N}$ the local (edge-weighted) mass operator 
\begin{equation}
	\label{eq:mass}
	M_x\colon \hat V \to \hat V,\qquad (M_x u, v) \coloneqq \frac{1}{2}\sum_{y \sim x} |x - y|^{2-\alpha}\, u(x)v(x),
\end{equation}
where $(\cdot,\cdot)$ denotes the Euclidean inner product. 
The diagonal global mass operator $M\colon \hat V \to \hat V$ is then defined as the sum of the local contributions $M\coloneqq \sum_{x\in\mathcal{N}}M_x$.
The bilinear form $\vspf{M \cdot}{\cdot}$ is an inner product on the space $\hat V$ with induced norm $\mnormf{\cdot}^2 \coloneqq \vsp{M \cdot}{\cdot}$. 
For subsets $S \subset \mathcal N$, we define localized versions of the mass operator by
\begin{equation}
	\label{eq:locM}
	M_S \coloneqq \sum_{x \in S} M_x
\end{equation}
and define the norm $\mnormfo{\cdot}{S}^2 \coloneqq \vsp{M_S \cdot}{\cdot}$. The local (edge-weighted) stiffness operator is given by 
\begin{equation}
	\label{eq:wgraphlaplacian}
	L_x\colon \hat V \to \hat V,\qquad  (L_x u, v) \coloneqq  \frac{1}{2}\sum_{y\sim x} \frac{[u(x) - u(y)][v(x) - v(y)]}{|x - y|^\alpha}.
\end{equation}
and the global stiffness operator is obtained by  $L\coloneqq \sum_{x\in\mathcal{N}}L_x$. 
For subsets $S \subset \mathcal{N}$, we define localized versions $L_S$ of the operator $L$ analogously to~\cref{eq:locM}. Note that $(L_S v)(x)$ is non-zero for nodes $x$ outside of $S$ that are adjacent to nodes in $S$. 
Since~$\mathcal G$ is connected, the kernel of $L$ consists only of the globally constant functions. Combined with the fact that there is at least one Dirichlet node, the bilinear form $(L\cdot,\cdot)$ defines an inner product on~$V$ with induced norm $\lnormf{\cdot}^2 \coloneqq (L\cdot, \cdot)$. Localized versions of this norm for $S \subset \mathcal{N}$ are denoted by $|\cdot|_{L,S}^2\coloneqq (L_S \cdot,\cdot)$.

Given a source term $f \in \hat V$, our model problem seeks a function $u \in V$ that satisfies
\begin{equation}\label{eq:weakform}
	(K u, v) = (f, v)\qquad \text{for all }v \in V.
\end{equation}
Here, the operator $K \colon \hat V \to \hat V$ is an abstract operator with the following properties.

\begin{samepage}
	\begin{assumption}[Properties of $K$]\label{ass:K}
	The operator $K$ satisfies the following conditions.
	\begin{enumerate}
		\item \emph{(Spectral bounds):}
		\label{item:speq} There exists constants $0<\gamma \leq \gamma^\prime<\infty$ such that
		 \begin{equation*}
			\gamma\, (Lv,v)\leq (Kv,v) \leq \gamma^\prime\, (Lv,v)\qquad\text{for all } v \in V.
		\end{equation*}
		\item \emph{(Sparsity pattern):}
		\label{item:locdec} $K$ only includes communication between neighboring nodes. In particular, it shares the sparsity pattern of $L$ and can be decomposed as
		\begin{equation*}
			K = \sum_{x \in \mathcal N} K_x,
		\end{equation*}
		where $K_x\colon \hat V \to \hat V$ are symmetric positive semi-definite operators with support only on $x$ and its neighboring nodes.
	\end{enumerate}
\end{assumption}  
\end{samepage}

Note that it is an immediate consequence of \cref{ass:K}, \cref{item:locdec} that $K$ is symmetric. By \cref{ass:K}, \cref{item:speq} it holds that $(K\cdot,\cdot)$ is an inner product on $V$ with induced norm $|\cdot|_K^2 \coloneqq (K\cdot,\cdot)$. For subsets $S \subset \mathcal N$ we further define the localized counterpart $K_S$ of $K$ analogously to \cref{eq:locM}.
The existence of a unique solution to this problem can be concluded from the fact that $(K\cdot,\cdot)$ is an inner product on $V$. Noting that $(f,\cdot) \in V^\prime$, where  $V^\prime$ denotes the dual space of $V$, the unique solvability follows from the Riesz representation theorem. Since $\mathcal N \cap \Gamma \neq \emptyset$, there exists a Friedrichs constant $C_\mathrm{fr}>0$ such that $|v|_M \leq C_\mathrm{fr} |v|_L$ for all $v \in V$. Denoting $|f|_{M^{-1}}\coloneqq \sup_{|v|_M = 1}(f,v)$, an immediate consequence of this inequality is the following  stability estimate for the solution~$u$ to~\cref{eq:weakform}, 
\begin{equation*}
	\lnormf{u} \leq C_\mathrm{fr}\gamma^{-1} |f|_{M^{-1}}.
\end{equation*}
Note that the dependence of the Friedrichs constant on the considered spatial network is generally unknown. However, we have the characterization $C_\mathrm{fr} = \lambda_1^{-1/2}$, where $\lambda_1$ is the smallest eigenvalue of the graph Laplacian problem $Lv = \lambda Mv$ subject to homogeneous Dirichlet boundary conditions. Under locality assumptions on the spatial network, one can relate $\lambda_1$ to the first eigenvalue of the normalized graph Laplacian, which has been extensively studied in the literature; see, e.g.,~\cite{CY95,C05}. 

\subsection{Possible choices of $K$}\label{sec:choiceK}

This subsection will provide an insight into exemplary classes of problems which can be considered with our formulation.
\begin{enumerate}
	\item \emph{(Finite differences and finite elements)}.\; Our formulation includes finite difference and finite element discretizations of second order linear elliptic operators. The finite difference case is covered for the choice $\alpha = 2$, while the finite element case is covered for shape-regular meshes and the choice $\alpha = 2-n$, where $n$ denotes the spatial dimension. For the case $n \geq 3$, the choice $\alpha = 2-n$ would violate the condition $\alpha \in [0,2]$. Therefore, we need the additional assumption that the considered hierarchy of meshes is quasi-uniform, which allows us to  rescale the equation with the minimal mesh size.
	\item \emph{(Weighted graph Laplacian)}.\; Also included in our formulation is the weighted graph Laplacian on the spatial network $\mathcal G$. Given $\alpha \in [0,2]$ and uniformly positive and  bounded edge-weights $\gamma_{xy}$, the operator $K$ in this case reads	
	\begin{equation*}
	 (Ku, v) \coloneqq \sum_{x \in \mathcal N}  \frac{1}{2}\sum_{y\sim x}\gamma_{xy} \frac{[u(x) - u(y)][v(x) - v(y)]}{|x - y|^\alpha}.
	\end{equation*}
	For the choice $\alpha = 1$, this model can be interpreted as a finite element discretization of a diffusion-type problem on the spatial network, involving one-dimensional problems at the edges and appropriate continuity conditions for the solution and its flux at the nodes; see also \cite{hellman2022wellposedness}. 
\end{enumerate}

We emphasize that also many other problems, such as spring or beam models of a fibre network, can be considered with our formulation; see, e.g., \cite{EGHKM24,KMM20}. However, these problems are not considered in the numerical examples in this paper and we therefore forego a detailed discussion.

\section{Coarse partition and assumptions on the network}\label{sec:network}

In this section, we first introduce some basic notation for graphs as well as (coarse) partitions of the graph. Further, we pose necessary assumptions on the network. We also provide some tools that will be essential for our construction of an appropriate coarse-scale approximation of the solution~$u$ to~\cref{eq:weakform} in the next section.  

For a subset of nodes $S \subset \mathcal N$, we denote with $\operatorname{dist}_S(x,y)$ the graph distance of the nodes $x,y\in S$ with respect to $S$, which corresponds to the smallest sum of edge lengths among all paths consisting of edges with nodes in $S$ connecting $x$ and $y$. If no such patch exists, we set $\operatorname{dist}_S(x,y) \coloneqq \infty$. Given another subset $R \subset S$,
\begin{equation}
	\label{eq:dist}
	\operatorname{dist}_S(x,R) = \min_{y \in R} \operatorname{dist}_S(x,y)
\end{equation}
denotes the graph distance between $x$ and $R$ considering only edges with both nodes in~$S$.

\subsection{Partitions and patches}

The definition of the proposed method is based on a partition of the set of nodes $\mathcal N$ into pairwise disjoint subsets. Mimicking the notation of finite element meshes, we denote the set of subsets by $\TH$ and might call the subsets \emph{elements}. We assume that for any partition, the corresponding subgraphs are connected. Note that this may be a direct consequence of the partitioning algorithm used; see~\cref{sec:numexp} below. Once again in analogy to finite element meshes, the parameter $H>0$ indicates the maximal diameter of the elements with respect to the graph distance, i.e.,
\begin{equation}
	\label{eq:H}
	H \coloneqq \max_{T \in \TH} H_T\quad\text{with}\quad H_T \coloneqq \max_{x,y\in T}\operatorname{dist}_T(x,y).
\end{equation}
Note that we use the parameter $H$ instead of $h$ to indicate the coarseness of the partition. 
Next, we introduce the concept of patches. Given a union of elements~$S$, we define the corresponding \emph{first-order patch} $\mathsf N(S)$ by
\begin{equation*}
	\mathsf N (S)\coloneqq  \bigcup \{T \in \TH \with T \subset S\;\;  \text{or}\;\; \exists e=\{x,y\}\in \mathcal E\with x \in T,\;y \in S\}.
\end{equation*}
For $\ell\in \mathbb N$, the \emph{$\ell$-th order patch} of $S$ is then recursively given by
\begin{equation*}
	\mathsf N^1(S) \coloneqq \mathsf N(S),\qquad \mathsf N^\ell(S)\coloneqq \mathsf N(\mathsf N^{\ell-1}(S)),\quad \ell \geq 2.
\end{equation*}

\subsection{Network assumptions}
For the analysis below, we need a hierarchy of partitions $\{\TH\}_{H \in \mathcal H}$ for some finite set of parameters $\mathcal H$ with minimal parameter $H_\mathrm{min}>0$. 
In particular, we require that the networks of interest behave similarly to continuous objects for sufficiently coarse $H_\mathrm{min}$. For example, we need an element-wise Poincar\'{e} inequality. To make this rigorous, we will make two assumptions on the hierarchy of partitions. 

The first one concerns an \emph{isoperimetric inequality} on the partitions, which involves the so-called \emph{isoperimetric dimension}. The isoperimetric dimension of a graph, denoted by $d\in \mathbb N$, describes how its behavior resembles that of an Euclidean space. Note that the isoperimetric dimension $d$ in general is not equal to the dimension $n$ of the domain~$\Omega$. To state the precise assumption, we need some additional notation. For $S \subset \mathcal N$, the degree of the node $x \in S$ is denoted by $\operatorname{deg}_S(x)$, which counts the number of nodes in $S$ that are adjacent to $x$. We then define the volume of $R \subset S \subset \mathcal N$ by
\begin{align*}
	\operatorname{vol}_S(R) \coloneqq \sum_{x \in R} \operatorname{deg}_S(x).
\end{align*}
Furthermore, we denote the boundary nodes of $R$ by
\begin{equation*}
	\operatorname{bnd}_S(R)\coloneqq \{x \in R\with \exists e = \{x,y\} \in \mathcal E \with y \in S\setminus R \}.
\end{equation*}
The definition of boundary nodes can be directly transferred to elements $T \in \TH$, where we also include nodes on the boundary of $\Omega$. That is,
\begin{equation*}
	\partial T \coloneqq \big(\partial \Omega \cap T\big) \bigcup \operatorname{bnd}_{\mathcal N}(T).
\end{equation*}
We can now formulate the isoperimetric inequality for the elements of the partitions.
 \begin{assumption}[Isoperimetric inequality]
	\label{ass:isoperimetric}
	There exist constants $\nu, \nu^\prime>0$ independent of $H$ such that
	\begin{equation*}
		(\operatorname{vol}_T(T))^{1/d}\leq \nu \frac{H}{H_\mathrm{min}}
	\end{equation*}
	for all $T \in \TH$,
	and the isoperimetric inequality 
	\begin{equation*}
		(\operatorname{vol}(R))^{(d-1)/d}\leq \nu^\prime\#(\operatorname{bnd}_T(R))
	\end{equation*}
	holds for all $R \subset T \in \TH$ with $\operatorname{vol}(R)\leq \operatorname{vol}(T\setminus R)$.
\end{assumption}
The isoperimetric constant is a measure of how well the nodes within the elements are connected.
Note that since we are considering finite graphs, an isoperimetric inequality is satisfied for any isoperimetric dimension. However, an inappropriate choice of the isoperimetric dimension may lead to very large constants.
If a graph has isoperimetric dimension $d$ it essentially means that the graph mimics a $d$-dimensional manifold. In the numerical experiments in \cref{sec:numexp} we consider graphs of isoperimetric dimension $d = 2.$

The second assumption on the network concerns the shape of the partitions and, once again, has connections to classical assumptions in the context of finite elements. 

\begin{assumption}[Appropriate partitions]
	\label{ass:network}
	There exist constants $\sigma,\rho, \rho^\prime >0$ independent of $H$ such that the following inequalities hold.
	\begin{enumerate} 
		\item \emph{(Quasi-uniformity):} \label{item:qu} For any $T \in \TH$, we have that
		\begin{equation*}
			\max_{T \in \TH} |1|_{M,T}^2 \leq \sigma \min_{T \in \TH} |1|_{M,T}^2.
		\end{equation*} 
		\item \emph{(Non-degenerateness):} \label{ass:nd} For any $T \in \TH$ and any subset $\hat T \subset T$, it holds that
		\begin{equation*}
			\rho\,|1|_{M,T}^2 \leq  |1|_{M,\hat T}^2\quad\text{and}\quad \operatorname{dist}(x,\partial T) \geq \rho^\prime H \quad\text{for all }  x \in \hat T.  
		\end{equation*}
		\item \emph{(Locality):}
		\label{ass:locality} We have that
		\begin{equation*}
			\max_{\{x,y\}\in \mathcal E} |x-y|\leq H_\mathrm{min}.
		\end{equation*}
	\end{enumerate}
\end{assumption}
The above assumptions are very similar to their counterparts in finite element theory. More specifically, \cref{item:qu} states that all elements are of similar size and \cref{ass:nd} prevents the elements from degenerating. \cref{ass:locality} is a bound on the maximal edge length, which can be easily enforced by refining certain edges and adding additional nodes. 

\subsection{Element-wise Poincar\'{e} inequality}

A very important property for the subsequent analysis is an element-wise Poincar\'e inequality. With the element average functional
\begin{equation}
	\label{eq:elemav}
	q_T(v)\coloneqq \frac{\vsp{M_T v}{1}}{\mnormfo{1}{T}^2},\qquad v \in \hat V,
\end{equation}
we can prove the following result. 
	
	\begin{samepage}
		\begin{lemma}[Element-wise Poincar\'{e} inequality]
			\label{l:poincare}
			For all $T \in \TH$ there exists a constant $C_\mathrm{po}>0$ such that 
			\begin{equation*}
				|v - q_T(v)|_{M,T} \le C_\mathrm{po} |v|_{L,T}
			\end{equation*}
			for all $v \in \hat V$. Provided that~\Cref{ass:network,ass:isoperimetric} are satisfied, there exists a constant $\mu>0$ independent of $H$ such that
			\begin{equation*}
				C_\mathrm{po} \leq \mu H.
			\end{equation*}
		\end{lemma}
	\end{samepage}

\begin{proof}
	A proof of this result can, e.g., be found in~\cite[Lem.~3.5 \& 3.6]{GoHeMa22}. 
\end{proof}

Note that for any reasonable partition, the above assumptions hold for a finite difference or finite element discretization. The corresponding grids or meshes are  well connected, so \cref{ass:isoperimetric} holds with a small constant. Uniform refinement produces grids or meshes with the same properties, and thus the isoperimetric dimension is preserved, and \cref{ass:isoperimetric} holds for the refined mesh with a similar constant. 
Similarly, given a reasonable partition, the above assumptions also hold for a large class of randomly generated spatial networks; see \cite{GoHeMa22,HaM22} for a numerical study of the $H$-dependence of the Poincar\'{e} constant for the case $\alpha = 1$. 
		
\subsection{Projection operators}
An important ingredient for the coarse-scale approximation introduced below is the $(\cdot,\cdot)_M$-orthogonal projection operator $\PiH \colon \hat V \to \Pnull(\TH)$ onto the space of $\TH$-piecewise constants defined by
\begin{equation*}
	\Pnull(\T_H)\coloneqq \operatorname{span}\{\one_\elem\with \elem\in \TH\} \subset \hat V,
\end{equation*}
where $\one_\elem$ denotes the indicator functions of the element $T \in \TH$. Recalling  \cref{eq:elemav}, the $L^2$-type projection $\PiH$ can be characterized as follows
\begin{equation}\label{eq:defL2}
	\PiH v =  \sum_{\elem \in \mathcal T_H}q_T(v)\one_\elem.
\end{equation}
We also introduce a bubble operator which, unlike $\PiH$, satisfies a stability estimate in the $L$-norm. It is constructed similarly to~\cref{eq:defL2} where the indicator functions $\{\one_\elem\}_{T \in \TH}$ are replaced by element-local bubble functions $\{b_T\}_{T \in \TH}$ that have the same element averages, cf.~\cref{eq:elemav}. We define the bubble function $b_T\in \V$ corresponding to the element $T \in \TH$ by
\begin{equation}
	\label{eq:bubble}
	b_T \coloneqq \frac{\tilde b_T}{q_T(\tilde b_T)},\qquad \tilde b_T(x) = \left\{\begin{array}{ll}
		\operatorname{dist}_T(x,\partial T)& \text{if }x \in T,\\
		0 & \text{else}.
	\end{array}\right.
\end{equation}
The definition of the bubble operator $\mathcal B_H\colon  \hat V \to \operatorname{span}\{b_T\with T \in \TH\}$ now reads
\begin{equation}
	\label{eq:bubbleop}
	\mathcal B_H v \coloneqq  \sum_{\elem \in \mathcal T_H}q_T(v)b_\elem.
\end{equation}
Note that, by definition, $\operatorname{ker} \PiH = \operatorname{ker} \mathcal B_H$. 
The following lemma states local stability and approximation results for $\PiH$ and~$\mathcal B_H$.  
\begin{lemma}[Properties of $\PiH$ and~$\mathcal B_H$]\label{l:projop}
	Let~\Cref{ass:network,ass:isoperimetric} be satisfied.
	Then, there exists a constant $C_\Pi>0$ independent of $H$ such that for any $T \in \TH$ and all $v \in \hat V$
	\begin{align}
		\label{eq:L2}
		\begin{split}
			|\PiH v|_{M,T}&\leq |v|_{M,T},\\
			|v-\PiH v|_{M,T}&\leq C_\Pi H |v|_{L,T}.
		\end{split}
	\end{align}
	Furthermore, there exists a constant $C_{\mathcal B}>0$ independent of $H$ such that for all $T \in \TH$ and $v \in \hat V$
	\begin{align}
		\label{eq:bubbleest}
		|\mathcal B_H v|_{M,T} + H|\mathcal B_H v|_{L,T} \leq C_{\mathcal B} |v|_{M,T}.
	\end{align}
\end{lemma}

\begin{proof}
	The first estimate in~\cref{eq:L2} follows directly from the orthogonality of the $L^2$-type projection $\Pi_H$. The second estimate follows from~\cref{l:poincare} with $C_\Pi = \mu$.

	To prove~\cref{eq:bubbleest}, we first derive estimates for the non-scaled bubble function $\tilde b_T$ defined in~\cref{eq:bubble}. Recalling \cref{eq:H}, we directly get that 
	\begin{equation}\label{eq:boundTbT}
		|\tilde b_T|_{M,T}\leq H|1|_{M,T}.
	\end{equation}
	With \cref{eq:boundTbT,ass:network}, we compute for some $\hat T \subset T$
	\begin{equation}
		\label{eq:esttildeb}
		\rho\rho^\prime H|1|_{M,T}^2\leq \rho^\prime H|1|_{M,\hat T}^2\leq (M_T\tilde b_T,1) \leq |\tilde b_T|_{M,T}|1|_{M,T} \leq H |1|_{M,T}^2.
	\end{equation}
	To estimate the $L$-norm of $\tilde b_T$, we note that $|\tilde b_T(x)-\tilde b_T(y)|\leq |x-y|$ for any edge $\{x,y\}\in \mathcal E$, which is a consequence of the triangle inequality. Therefore, we obtain
	\begin{equation}
		\label{eq:lnormtb}
		|\tilde b_T|_{L,T}^2 = \frac12\sum_{x \in T}\sum_{y \sim x}\frac{[\tilde b_T(x)-\tilde b_T(y)]^2}{|x-y|^\alpha}\leq  \frac12\sum_{x \in T}\sum_{y \sim x} |x-y|^{2-\alpha} = |1|_{M,T}^2.
	\end{equation}
	Using \cref{eq:esttildeb,eq:elemav}, and once again~\cref{eq:boundTbT}, we get that
	\begin{align*}
		|\mathcal B_H v|_{M,T} = \frac{|q_T(v)|}{|q_T(\tilde b_T)|}|\tilde b_T|_{M,T} \leq (\rho\rho^\prime)^{-1}|v|_{M,T}.
	\end{align*}
	Similarly, we obtain using~\cref{eq:esttildeb,eq:elemav,eq:lnormtb} that 
	\begin{align*}
		|\mathcal B_H v|_{L,T} = \frac{|q_T(v)|}{|q_T(\tilde b_T)|}|\tilde b_T|_{L,T} \leq 	(\rho\rho^\prime)^{-1} H^{-1}|v|_{M,T},
	\end{align*}
	such that~\cref{eq:bubbleest} follows with the constant $C_{\mathcal B} = (\rho\rho^\prime)^{-1}$.
\end{proof}

\section{Prototypical coarse-scale approximation}\label{sec:protomethod}

In this section, we introduce a prototypical approach for a reliable coarse-scale approximation of the solution to~\cref{eq:weakform}. We use the term \emph{prototypical} to indicate that the method is not intended for practical use, but only for an initial theoretical study. A practical method is presented in the two subsequent sections.

The construction of the prototypical method involves the notion of \emph{corrections}. 
We introduce the space of fine-scale functions $\W \coloneqq \operatorname{ker}\PiH|_V$ which consists of functions that are not captured by the coarse-scale operator $\PiH$. The correction operator $\C \colon V \to \W$ then assigns to any $v \in V$ the unique function $\C v \in \W$ that solves
\begin{equation}
	(K\C v,\,w) = (Kv,\,w)\qquad \text{for all }w\in \W.
\end{equation}
We now construct the approximation space for the prototypical method by adding fine-scale corrections to bubble functions, i.e., we define
\begin{equation}
	\label{eq:idealsp}
	V_H \coloneqq (1-\C)\mathcal B_H \Pnull(\TH).
\end{equation}
Note that $\dim V_H = \#\TH$ and the functions $\{(1-\C)b_T\}_{T \in \TH}$ are a basis of $V_H$. In the following, the space~$V_H$ is used as a problem-adapted approximation space. 
More precisely, the prototypical method is defined as a Galerkin approach in the space $V_H$, i.e., we seek $u_H\in V_H$ such that 
\begin{equation}\label{e:galerkinideal}
	(K u_H,\,v) = ( f,\,v)\qquad\text{for all }v_H \in V_H.
\end{equation}
The following approximation result provides a justification of the above construction.

\begin{lemma}[Uniform  approximation]\label{l:ua}
	Let \Cref{ass:network,ass:isoperimetric} be satisfied. Then, for any $f \in \hat V$, the  prototypical approximation $u_H$ defined in \eqref{e:galerkinideal} satisfies
	\begin{equation}
		\label{eq:firstorderestideal}
		\lnormf{u-u_H} 
		\leq C_\Pi \gamma^{-1}H |f|_{M^{-1}}.
	\end{equation}
	Denoting  $\tilde f \coloneqq M^{-1}f$ and recalling that
	\begin{equation*}
		|f|_{M^{-1}}= \sup_{|v|_M = 1}(f,v) = |\tilde f|_M,
	\end{equation*}
	we also obtain the second-order error estimate
	\begin{equation}
		\label{eq:estidealmethod}
		\lnormf{u-u_H} 
		\leq C_\Pi \gamma^{-1}H \mnormf{\tilde f-\PiH \tilde f}
		\leq C_\Pi^2 \gamma^{-1}H^2\lnormf{\tilde f}.
	\end{equation}
\end{lemma}

\begin{proof}
	Since $e \coloneqq u-u_H \in \W$ by the Galerkin orthogonality, we obtain using \cref{ass:K} and \cref{l:projop} that 
	\begin{align*}
		\gamma\lnormf{e}^2  &\le \vsp{Ke}{e} = \vsp{Ku}{e} = (f,e) \leq |f|_{M^{-1}}|e|_M
		\leq C_\Pi H |f|_{M^{-1}}|e|_L,
	\end{align*}
	which proves \cref{eq:firstorderestideal}. For the proof of the second-order estimate \cref{eq:estidealmethod}, we employ the orthogonality properties of $\PiH$ and obtain 
	\begin{align*}
		(f,e) &= (M\tilde f,e-\PiH e) 
		= \vspf{M(\tilde f-\PiH \tilde f)}{e - \PiH e}\leq \mnormf{\tilde f-\PiH \tilde f}\mnormf{e - \PiH e}\\
		&
		\leq C_\Pi H\mnormf{\tilde f-\PiH \tilde f} |e|_L \leq C_\Pi^2H^2|\tilde f|_L|e|_L,
	\end{align*}
	where we once again applied~\cref{l:projop}. This proves~\cref{eq:estidealmethod}.
\end{proof}

\begin{remark}[$\TH$-piecewise source terms]
	\label{rem:exactnessideal}
	For $\tilde f\in \Pnull(\TH)$, the prototypical method is exact, as for $\TH$-piecewise constants $\mnormf{\tilde f-\PiH \tilde f} = 0$ in \cref{eq:estidealmethod}.
\end{remark}

\section{Exponential decay and localization}\label{sec:decloc}
In this section, we prove the exponential decay behavior of the corrections, which motivates a localization of the basis functions. This is used in the following section to construct a practical localized version of the prototypical method introduced above.

A first observation is that the correction operator $\C$ can be split into element-wise contributions, that is
\begin{equation}
	\label{eq:sumelemcorr}
	\C = \sum_{\elem \in \mathcal \TH} \C_T,
\end{equation} 
where the operator $\C_T\colon V\to \W$ for any $T \in \TH$ assigns to a function $v \in V$ the unique solution $\C_Tv \in \W$ to
\begin{equation}\label{def:phi}
	(K\C_T v,w) = (K_T v,w)\qquad \text{for all }w \in \W.
\end{equation}
Recall that the $K_T$ is the localization of the operator $K$ to the element $T$, which is defined analogously to~\cref{eq:locM}. 
In analogy to a finite element-based discretization in the continuous setting, we call the operators $\{\C_T\}_{T \in \TH}$ \emph{element corrections}. The following theorem states an exponential decay of the element corrections away from their associated elements.

\begin{theorem}[Exponential decay]
	\label{thm:decay}
	Let \Cref{ass:network,ass:isoperimetric} be satisfied. Then there exists a constant $c>0$ independent of $H$ such that for any $T \in \TH$ and all $\ell \in \mathbb N$
	\begin{equation*}
		|\C_T v|_{L,\mathcal N \setminus \mathsf N^\ell(T)} \leq \exp(-c\ell)|\C_T v|_L\quad\text{for all } v \in V.
	\end{equation*}
\end{theorem}
\begin{proof}
	We abbreviate $\varphi \coloneqq \C_T v$. 
	With \Cref{ass:network,ass:isoperimetric}, one can prove that a path of length $H$ in graph distance crosses at most a fixed number of partitions $k\in \mathbb N$, where $k$ depends only on the constants $\nu,\nu^\prime,\sigma,\rho,$ and $\rho^\prime$. 
	
	Similarly as in the continuous setting (see, e.g., \cite[Thm.~3.15]{Peterseim2021}), the proof employs cut-off techniques. Consider $\ell \geq k$ and denote by
	\begin{equation*} 
		\Gamma_1 \coloneqq \operatorname{bnd}_{\mathsf N^\ell(T)}(R),\qquad \Gamma_2 \coloneqq \operatorname{bnd}_{\mathcal N \setminus \mathsf N^{\ell-k}(T)}(R)
	\end{equation*}
	the inner and outer boundary of the layer $R \coloneqq \mathsf N^{\ell}(T) \setminus \mathsf N^{\ell-k}(T)$, respectively. We then define the cut-off function $\eta \in \hat V$ by setting
	\begin{equation*}
		\eta(x) \coloneqq \left\{\begin{array}{ll}
			\frac{\operatorname{dist}(x,\Gamma_2)}{\operatorname{dist}(x,\Gamma_1) + \operatorname{dist}(x,\Gamma_2)}& \text{if }x \in R,\\
			0 & \text{if }x \in \mathsf N^{\ell-k}(T),\\
			1 & \text{else}.
		\end{array}\right.
	\end{equation*}
	Due to the choice of $k$, there exists a constant $C_\eta>0$ such that
	\begin{equation}
		\label{eq:eta}
		|\eta(x)-\eta(y)| \leq C_\eta H^{-1} |x-y|
	\end{equation}
	for any edge $\{x,y\} \in \mathcal E$.
	
	Using the spectral equivalence of $K$ and $L$ (cf.~\cref{ass:K}, \cref{item:speq}) and recalling that $\varphi \in \W$, we obtain by the construction of $\eta$ that
	\begin{align*}
		\gamma |\varphi|_{L,\mathcal N \setminus \mathsf N^\ell(T)}^2 = \gamma |(1-\mathcal B_H)\varphi|_{L,\mathcal N \setminus \mathsf N^\ell(T)}^2 \leq \gamma |(1-\mathcal B_H)(\eta\varphi)|_{L}^2 \leq  |(1-\mathcal B_H)(\eta\varphi)|_{K}^2.
	\end{align*}
	Abbreviating $w \coloneqq (1-\mathcal B_H)(\eta\varphi) \in \W$, we get
	\begin{align*}
		  |(1-\mathcal B_H)(\eta\varphi)|_{K}^2 &= (K(1-\mathcal B_H)(\eta\varphi),w)\\
		  &= (K(1-\mathcal B_H)\varphi,w)-(K(1-\mathcal B_H)(1-\eta)\varphi,w) \eqqcolon \Xi_1 + \Xi_2.
	\end{align*}
	We first consider the term $\Xi_1$. Using  that $\varphi\in \W$ and  $\operatorname{supp}(w) \subset \mathcal N \setminus \mathsf N^{\ell-k}(T)$ by the definition of $\eta$ and the locality of the operator $\mathcal B_H$, we obtain with~\eqref{def:phi}
	\begin{align*}
		\Xi_1 = (K\varphi,w) = (K_T v,w) = 0.
	\end{align*} 
	To estimate the term $\Xi_2$, we use the spectral equivalence of $K$ and $L$ (cf.~\cref{ass:K}, \cref{item:speq}) and the locality of the operator $K$ (cf.\cref{ass:K}, \cref{item:locdec}). By the choice of~$\eta$, we get
	\begin{align}\label{eq:theta2}
		|\Xi_2| \leq \gamma^\prime |(1-\mathcal B_H)(1-\eta)\varphi|_{L,R}|(1-\mathcal B_H)(\eta\varphi)|_{L,R}.
	\end{align}
	Both factors on the right-hand side can be estimated similarly. Therefore, we only consider the second factor. Using the triangle inequality, we obtain 
	\begin{align}
		\label{eq:proofterms}
		|(1-\mathcal B_H)(\eta\varphi)|_{L,R}\leq |\eta\varphi|_{L,R} + |\mathcal B_H (\eta\varphi)|_{L,R}
	\end{align}
	To estimate $|\eta\varphi|_{L,R}$ we denote by $\bar R$ the set of nodes which are in $R$ or neighbored to a node in $R$. Using a discrete analogon of the product rule in the continuous setting, \cref{eq:eta}, and the discrete Cauchy--Schwarz inequality, we get that
	\begin{align*}
		|\eta\varphi|_{L,R}^2 &\leq \sum_{x \in \bar R}\frac12 \sum_{\bar R \ni y \sim x}\frac{((\eta\varphi)(x)-(\eta\varphi)(y))^2}{|x-y|^\alpha} \\
		&= \sum_{x \in \bar R}\frac12 \sum_{\bar R \ni y \sim x}\eta(x)\frac{[(\eta\varphi)(x)-(\eta\varphi)(y)][\varphi(x)-\varphi(y)]}{|x-y|^\alpha}\\
		&\qquad  + \sum_{x \in \bar R}\frac12 \sum_{\bar R \ni y \sim x} \varphi(x)\frac{[(\eta\varphi)(x)-(\eta\varphi)(y)][\eta(x)-\eta(y)]}{|x-y|^\alpha}\\
		&\leq 2|\eta\varphi|_{L,R}|\varphi|_{L,R} + 2C_\eta H^{-1}|\varphi|_{M,R}|\eta\varphi|_{L,R}.
	\end{align*}
	Note that in order to swap $x$ and $y$ in the second sum, we have to add additional summands so that we sum over each edge twice. The factor of two is necessary to get back to the original definitions of $|\cdot|_{M,R}$ and $|\cdot|_{L,R}$. Dividing by $|\varphi|_{L,R}$ and using the fact that~$\varphi \in \W$ as well as \cref{eq:eta,l:projop}, we can estimate the first term on the right-hand side of~\cref{eq:proofterms} by
	\begin{equation}\label{eq:decay1}
		|\eta\varphi|_{L,R}\leq 2(1+C_\eta C_\Pi)|\varphi|_{L,R}.
	\end{equation}
	To estimate the second term on the right-hand side of \cref{eq:proofterms}, we again use that $\varphi \in \W$ and \cref{l:projop} to get that 
	\begin{align}\label{eq:decay2}
		|\mathcal B_H(\eta\varphi)|_{L,R}\leq C_{\mathcal B} H^{-1}|\eta \varphi|_{M,R}\leq C_{\mathcal B}C_\Pi |\varphi|_{L,R}.
	\end{align}
	We go back to~\eqref{eq:proofterms} and combine~\cref{eq:decay1,eq:decay2}, which yields
	\begin{equation*}
		|(1-\mathcal B_H)(\eta\varphi)|_{L,R} \leq (2(1+C_\eta C_\Pi) + C_{\mathcal B}C_\Pi)|\varphi|_{L,R}.
	\end{equation*}
	Finally, we obtain using the previous estimates (that can verbatim be applied to the second factor on the right-hand side of~\eqref{eq:theta2} as well)
	\begin{equation*}
		|\varphi|^2_{L,\mathcal N\setminus \mathsf{N}^{\ell}(T) } \leq  C 	|\varphi|^2_{L,\mathsf{N}^{\ell}(T) \setminus \mathsf{N}^{\ell-k}(T)} = C |\varphi|^2_{L,\mathcal N \setminus  \mathsf{N}^{\ell-k}(T)} -  C |\varphi|^2_{L,\mathcal N \setminus \mathsf{N}^{\ell}(T)},
	\end{equation*}
	where $C \coloneqq  \gamma^\prime (2(1+C_\eta C_\Pi) + C_{\mathcal B}C_\Pi)^2 /\gamma$. Therefore, we conclude that
	\begin{equation*}
		|\varphi|_{L,\mathcal{N} \setminus \mathsf{N}^{\ell}(T) } \leq \delta |\varphi|_{L,\mathcal N \setminus \mathsf{N}^{\ell-k}(T) } \leq \delta^{\lfloor \ell/k\rfloor} |\varphi|_{L}
	\end{equation*}
	with $\delta \coloneqq (\tfrac{C}{1+C})^{1/2}<1$. Note that this estimate remains valid for $\ell <k$. The assertion follows with  $\delta^{\lfloor \ell/k\rfloor}\leq \delta^{-1}(\delta^{1/k})^\ell\leq \sqrt{2}\exp( -\tfrac{1}{k}\log(1/\delta)\ell)$.
\end{proof}

The exponential decay results motivates the localization of the element corrections. Therefore, we introduce for $T \in \TH$ the local subspace $\W_T^\ell$ defined by
\begin{equation*}
	\W_T^\ell\coloneqq \{w \in \W \with \operatorname{supp}(w)\subset \mathsf N^\ell(T)\}.
\end{equation*} 
We then uniquely define the localized element corrector $\C_T^\ell\colon \V \to \W_T^\ell$ by seeking for any $v \in \V$ the solution $\C_T^\ell v \in W_T^\ell$ to
\begin{equation*}
	(K\C_T^\ell v,w) = (K_Tv,w)\qquad \text{for all }w \in \W_T^\ell.
\end{equation*}
In analogy to \cref{eq:sumelemcorr}, we define the localized correction operator as the sum of the element-wise contributions, i.e.,
\begin{equation*}
	\C^\ell v\coloneqq \sum_{\elem \in \mathcal T_H} \C_T^\ell v.
\end{equation*}
For the localized operator, the following exponential approximation result can be shown.

\begin{theorem}[Exponential approximation]\label{thm:expapp}
	Let \Cref{ass:network,ass:isoperimetric} be satisfied. Then there exists a constant $C_a>0$ independent of $H$ such that for any $T \in \TH$ and all $\ell \in \mathbb N$ 
	\begin{equation*}
		|\C v -\C^\ell v|_L \leq C_a\ell^{d/2} \exp(-c\ell)|v|_{L},
	\end{equation*}
	where $c$ is the constant from \cref{thm:decay}.
\end{theorem}

\begin{proof}
	The proof of this theorem uses cut-off techniques similar to those used in the proof of \cref{thm:decay}. For the sake of brevity, the proof is omitted here and the reader is referred to \cite[Lem.~4.6 \& 4.7]{GoHeMa22}. Therein, a similar result is proved for partitions~$\TH$ constructed based on Cartesian meshes.
\end{proof}

\section{Practical localized coarse-scale approximation}\label{sec:prac}
In this section, we present two practical versions of the prototypical method introduced in \cref{sec:protomethod}. Based on the above localization techniques, we construct approximation spaces with localized basis functions that can be obtained by solving local problems only. The first practical version is conceptually simple and already works with the above assumptions on the network, while the second version is more sophisticated but can achieve a more stable approximation under certain additional assumptions.

\subsection{Naive version}
The most straightforward way to obtain a localized version of the ansatz space defined in~\cref{eq:idealsp} is to replace the global correction operator $\C$ by its localized counterpart $\C^\ell$, i.e., 
\begin{equation}
	\label{eq:locsimple}
	V_H^\ell \coloneqq (1-\C^\ell)\mathcal B_H \Pnull(\TH).
\end{equation}
The corresponding practical method then seeks $u_H^\ell \in V_H^\ell$ that solves
\begin{equation}\label{eq:naiveloc}
	(Ku_H^\ell,v) = (f,v)\qquad \text{for all }v \in V_{H}^\ell.
\end{equation}
In the following theorem, we state a convergence result for this straight-forward localized method.

\begin{theorem}[Error estimate of the naive method]
	\label{thm:convsimple}
	Let \Cref{ass:network,ass:isoperimetric} be satisfied. Then, there exists a constant $C>0$ independent of $H$ such that for any source term $f \in \hat V$
	\begin{align}\label{eq:convsimple}
		\begin{split}
					|u-u_H^\ell|_L &\leq C\big( H |f|_{M^{-1}} + \ell^{d/2}H^{-1}\exp(-c\ell)|f|_{M^{-1}}\big)\\
			& \leq C\big(H^2|\tilde f|_L + \ell^{d/2}H^{-1}\exp(-c\ell)|f|_{M^{-1}}\big)
		\end{split}
	\end{align}
	with $\tilde f \coloneqq M^{-1}f$.
\end{theorem}

\begin{proof}
	By standard argument, we have quasi-optimality of the Galerkin solution $u_H^\ell$ in the sense that for any $v \in V_H^\ell$
	\begin{equation*}
		|u-u_H^\ell|_L \leq \gamma^{-1}\gamma^\prime |u-v|_L.
	\end{equation*}
	We choose $v \coloneqq (1-\C^\ell)\mathcal B_H u \in V_H^\ell$, add and subtract the prototypical approximation $u_H\in V_H$ from \cref{e:galerkinideal}, and use that due to $(1-\mathcal B_H)u\in \W$
	\begin{equation*}
		u_H = (1-\C)u = (1-\C)\mathcal B_H u + (1-\C)(1-\mathcal B_H)u = (1-\C)\mathcal B_Hu.
	\end{equation*}
	Altogether, we obtain with the triangle inequality
	\begin{align}
		\label{eq:locestsimple}
		|u-u_H^\ell|_L \leq \gamma^{-1}\gamma^\prime\big(|u-u_H|_{L} + |(\C-\C^\ell)\mathcal B_H u|_L\big).
	\end{align}
	The first term on the right-hand side can be estimated with~\cref{l:ua} and the second term is bounded by a combination of \cref{l:projop,thm:expapp}. 
\end{proof}

The error estimate in \cref{thm:convsimple} proves the optimal order convergence of the method provided that the localization parameter~$\ell$ is increased logarithmically with the mesh size. However, for a fixed localization parameter, the error may increase with decreasing mesh size due to the factor $H^{-1}$ in the term that quantifies the localization error, cf.~\cref{eq:convsimple}. A numerical illustration of this effect can be found in \cref{sec:numexp}. We address this issue in the next subsection. 

\subsection{Stabilized version}
To eliminate the negative power of the mesh size in the localization error, a more sophisticated localization strategy is needed. More precisely, we introduce an appropriate quasi-interpolation operator and proceed in a similar way as in~\cite{HaPe21,DHM23}. 
However, the challenge in the algebraic setting is the construction of a feasible quasi-interpolation operator. Unlike in the setting of triangular or quadrilateral meshes, operators based on piecewise linear functions are not available in the algebraic setting. Instead, we employ an abstract partition of unity which is assumed to satisfy the following properties.

\begin{assumption}[Partition of unity]\label{ass:pu}
	Suppose there exists an overlapping cover $\{U_T\}_{T \in \TH}$ of the node set $\mathcal N$ with $U_T \supset T$ for all $T \in \TH$ and corresponding partition of unity functions $\{\Lambda_T\}_{T \in \TH}$ with $\operatorname{supp}(\Lambda_T) \subset U_T$ such that the following requirements are satisfied.
	%
	\begin{enumerate}
		\item \emph{(Partition of unity and non-negativity):}
		\label{ass:sum} The functions $\Lambda_T$ are non-negative and satisfy
		\begin{equation*}
			\sum_{\elem \in \mathcal T_H}\Lambda_T \equiv 1.
		\end{equation*}
		\item \emph{(Lipschitz bound):}
		\label{item:lipschitz} There exists a constant $C_\Lambda>0$ independent of $H$ such that
		\begin{equation*}
			\max_{\{x,y\} \in \mathcal E}\frac{|\Lambda_T(x)-\Lambda_T(y)|}{|x-y|} \leq C_\Lambda H^{-1}.
		\end{equation*}
		\item \emph{(Bounded supports):}
		\label{item:boundedsupp} There exists some $k \in \mathbb N$ such that for all $T \in \TH$
		\begin{equation*}
			U_T \subset \mathsf N^k(T).
		\end{equation*}
	\end{enumerate}
\end{assumption}

Given a partition of unity that fulfills the above assumptions, we can define a quasi-interpolation operator $I_H \colon \hat V \to V$ as follows. For any $v \in \hat V$, we set 
\begin{equation}
	\label{eq:IH}
	I_H v\coloneqq \sum_{\elem \in \mathcal K_H} q_T(v) \Lambda_T,
\end{equation}
where $q_T$ is defined in~\cref{eq:elemav} and $\mathcal K_H \coloneqq \{T \in \TH \with \Lambda_T \in V\}$ denotes the set of elements whose corresponding partition of unity functions are $V$-conforming. 

To prove stability and approximation properties of $I_H$, we distinguish between elements for which the operator $I_H$ locally preserves constants and elements for which this is not the case. The set of elements for which constants are locally preserved is given by
\begin{equation*} 
	\mathcal G_H\coloneqq \{T \in \TH \with \forall K \in \TH \setminus \mathcal K_H \with T \cap U_K = \emptyset\}.
\end{equation*}
For elements $T \in \mathcal G_H$, we need to assume a local Poincar\'e inequality on $\mathsf N^k(T)$, while for elements $T \in \TH \setminus \mathcal G_H$, we need a local Friedrichs inequality. Since it is generally not guaranteed for $T \in \TH \setminus \mathcal G_H$ that $\mathsf N^k(T)$ has a non-empty intersection with $\Gamma$, which is required for the Friedrichs inequality, we consider the larger patch $\mathsf N^{2k}(T)$ in this case. 
These requirements are made more precise in the following assumption.

\begin{samepage}
\begin{assumption}[Poincar\'e and Friedrichs inequalities on patches]\label{ass:poincarefriedrichs}
	There exists a constant $\vartheta>0$ independent of $H$ such that for any $T \in \mathcal G_H$
	\begin{equation*}
		|v-q_{\mathsf N^k(T)}(v)|_{M,\mathsf N^k(T)} \leq \vartheta H|v|_{L,\mathsf N^k(T)}\quad\text{for all } v \in V,
	\end{equation*}
	where $q_{\mathsf N^k(T)}$ is defined as in~\cref{eq:elemav} but with $\mathsf N^k(T)$ instead of $T$. 
	Furthermore, there exists a constant $\vartheta^\prime>0$ independent of $H$ such that for any $T \in \TH \setminus \mathcal G_H$
	\begin{equation*}
		|v|_{M,\mathsf N^{2k}(T)} \leq \vartheta^\prime H|v|_{L,\mathsf N^{2k}(T)}\quad\text{for all } v \in V.
	\end{equation*}
\end{assumption}
\end{samepage}

Note that if \cref{ass:isoperimetric} holds for $k$-th and $2k$-th order patches of elements instead of just elements, then these Poincar\'e and Friedrichs inequalities can be proved similarly to \cref{l:poincare}. In the case of the Friedrichs inequality, we also need the boundary nodes to be sufficiently dense in $\Gamma$ so that $\mathsf N^{2k}(T)$ for $T \in \TH \setminus \mathcal G_H$ contains at least one boundary node.  We emphasize that \cref{ass:isoperimetric} in general does not imply the corresponding assumption for patches of elements.

We are now ready to state the stability and approximation properties of $I_H$.
 
\begin{lemma}[Properties of $I_H$]
	\label{l:propIH}
	Let \Cref{ass:network,ass:isoperimetric,ass:pu,ass:poincarefriedrichs} be satisfied. Then, there exists a constant  $C_I>0$ independent of $H$ such that 
	\begin{equation}
		\label{eq:iH}
		H^{-1}|v-I_H v|_{M} + |I_H v|_{L} \leq C_I |v|_{L}\quad\text{for all } v \in V.
	\end{equation}
	Moreover, $I_H$ only extends the support of a functions by at most $k$ layers with $k$ from \cref{ass:pu}. 
	%
\end{lemma}
\begin{proof}
	Let $T \in \TH$ be an arbitrary but fixed element. Using that the number of sets~$U_K$ that cover $T$ is bounded by~$C_\mathrm{ol}k^d$ for some constant $C_\mathrm{ol}>0$ independent of $H$ as well as 
	\cref{ass:network} (\cref{item:qu}) and \cref{ass:pu} (\cref{item:lipschitz}), we obtain for any $c \in \mathbb R$ that
	\begin{align}
		\label{eq:IHloc}
		\begin{split}
			|I_H(v-c)|_{L,T} &= \bigg| \sum_{U_K\supset T}q_K(v-c)\Lambda_K\bigg|_{L,T} \leq 	\sum_{U_K\supset T} |v-c|_{M,K}|1|_{M,K}^{-1}|\Lambda_K|_{L,T}\\
			&\leq C_\mathrm{ol}C_\Lambda \sigma  k^{d/2} H^{-1}  |v-c|_{M,\mathsf N^k(T)}.
		\end{split}
	\end{align}
	 In the following we distinguish between two cases: Firstly, we consider the case that $T \in \mathcal G_H$ for which $c = I_Hc$ holds on~$T$, i.e., constants are locally preserved. 	 
	 Choosing $c \coloneqq q_{\mathsf N^k(T)}(v)$ and  using the Poincar\'e inequality from \cref{ass:poincarefriedrichs} and the triangle inequality as well as $|c|_{L,T} = 0$, we  obtain that
	\begin{equation*}
		\label{eq:poincarepatch}
		|I_Hv|_{L,T} \leq	|I_H(v-c)|_{L,T} \leq C_\mathrm{ol}C_\Lambda \sigma\vartheta  k^{d/2} |v|_{L,\mathsf N^k(T)}.
	\end{equation*}
	Secondly, we consider the case that $T \in \TH \setminus \mathcal G_H$. In this case, we choose $c = 0$ and use the Friedrichs inequality from \cref{ass:poincarefriedrichs}. Similarly as before, we obtain
	\begin{equation*}
		|I_Hv|_{L,T}
		\leq C_\mathrm{ol}C_\Lambda \sigma\vartheta^\prime k^{d/2} |v|_{L,\mathsf N^{2k}(T)}.
	\end{equation*}
	The stability estimate in \cref{eq:iH} follows by summing over all elements.
	
	To show the approximation property in~\cref{eq:iH}, we use similar arguments as in~\cref{eq:IHloc} which yields for any $c \in \mathbb R$ that
	\begin{equation*}
		|I_H(v-c)|_{M,T} \leq C_\mathrm{ol} \sigma  k^{d/2} |v-c|_{M,\mathsf N^k(T)}.
	\end{equation*}
	The desired approximation result can then be concluded by
	\begin{align*}
		|v-I_H v|_{M,T} \leq |v-c|_{M,T} + |I_H(v-c)|_{M,T},
	\end{align*}
	considering the two different cases as above and summing over all elements.
\end{proof}

Based on $I_H$, we define another quasi-interpolation operator $P_H\colon \hat V \to V$ by
\begin{equation}
	\label{eq:PH}
	P_H v = I_H v  + \mathcal B_H(v-I_H v)\qquad\text{for all } v \in \hat V.
\end{equation}
The operator $P_H$ has favorable properties compared to $\PiH$ and $B_H$ regarding stability in~$|\cdot|_L$, while having the same kernel as $\PiH$. This is made precise in the following lemma.

\begin{lemma}[Properties of $P_H$]
	\label{l:propPH}
	Let \Cref{ass:network,ass:isoperimetric,ass:pu,ass:poincarefriedrichs} be satisfied. It holds that
	\begin{equation}
		\label{eq:kerneleq}
		\operatorname{ker} \PiH= \operatorname{ker}P_H	
	\end{equation}
	and $P_H$ admits the representation
	\begin{equation}
		\label{eq:repPH}
		P_H v = \sum_{\elem \in \mathcal T_H}q_T(v)P_H(b_T)\quad\text{for all } v \in \hat V.
	\end{equation}
	Furthermore, there exists a constant $C_P>0$ independent of $H$ such that
	\begin{equation}
		\label{eq:stabPH}
		|P_H v|_L \leq C_P |v|_L\quad\text{for any } v \in V
	\end{equation}
	and $P_H$ only extends the support of a functions by at most $k$ layers with $k$ from \cref{ass:pu}.
	\end{lemma}

	\begin{proof}
	The representation \cref{eq:repPH} is a direct consequence of \cref{eq:bubbleop} and the identities $\mathcal B_H^2 = \mathcal B_H$ and $I_H\circ \mathcal B_H = I_H$. The identity~\cref{eq:kerneleq} immediately follows from~\cref{eq:repPH}.
	
	To prove the stability estimate~\cref{eq:stabPH}, we use \cref{eq:PH,l:propIH,l:projop} and get
	\begin{align*}
		|P_Hv|_L \leq C_I |v|_L + C_{\mathcal B} H^{-1}|v-I_Hv|_M \leq C_I(1+C_{\mathcal B})|v|_L,
	\end{align*}
	which is the desired result.
\end{proof}

An alternative to defining a practical localized space as in \cref{eq:locsimple} is to use the operator $P_H$ instead of $\mathcal B_H$. This leads to the slightly modified approximation space
\begin{equation*}
	\tilde V_H^\ell \coloneqq (1-\C^\ell)P_H \Pnull(\TH).
\end{equation*}
Similarly as above, the corresponding practical method seeks the Galerkin approximation $\tilde u_H^\ell \in \tilde V_H^\ell$ that solves
\begin{equation*}
	(K\tilde u_{H,\ell},v) = (f,v)\qquad \text{for all }v \in \tilde V_{H}^\ell.
\end{equation*}
The following theorem quantifies the improvements of the stabilized version compared to the naive one defined in~\eqref{eq:naiveloc} in the sense that the term $H^{-1}$ is avoided in the error estimate, cf.~\cref{thm:convsimple}.

\begin{theorem}[Error estimate of the stabilized method]\label{thm:convstab}
	Let \Cref{ass:network,ass:isoperimetric,ass:pu,ass:poincarefriedrichs} be satisfied. 
	Then, there exists a constant $C>0$ independent of $H$ such that for any source term $f \in \hat V$
	\begin{align*}
		|u-\tilde u_H^\ell|_L &\leq C( H |f|_{M^{-1}} + \ell^{d/2}\exp(-c\ell)|f|_{M^{-1}})\\
		& \leq C(H^2|f|_L + \ell^{d/2}\exp(-c\ell)|f|_{M^{-1}}).
	\end{align*}
\end{theorem}

\begin{proof}
	We follow the lines from the proof of \cref{thm:convsimple}. Analogously to \cref{eq:locestsimple}, we obtain the estimate
	\begin{align*}
		|u-u_H^\ell|_L \leq \gamma^{-1}\gamma^\prime\big(|u-u_H|_{L} + |(\C-\C^\ell)P_H u|_L\big).
	\end{align*}
	The first term on the right-hand side can be bounded with \cref{l:ua}. The second term can be estimated by combining \cref{l:propPH,thm:expapp}.
\end{proof}

\section{Implementation and numerical examples}
\label{sec:numexp}

In this section, we first discuss some implementation details regarding the construction of a partition of a graph, which will be used to define a corresponding partition of unity. Thereafter, several numerical examples are presented which support the theoretical results of this paper.

\subsection{Implementation details}\label{sec:impl}
\begin{samepage}
\subsubsection*{Construction of a partition}
For the construction of a partition, there is great flexibility in which algorithm to use. Many graph partitioners implemented in different toolboxes could be used. \end{samepage} Note that typically \cref{ass:network} is not guaranteed a priori, but can be verified a posteriori.  For our implementation, we used the Gonzalez algorithm for the $N$-center clustering algorithm from \cite{Gon85,HoS85}, where $N = \#\TH$. The Gonzalez algorithm provides an approximate solution to the NP-hard $N$-center clustering problem, which seeks a set of centers $C = \{x_1,\dots,x_m\}$ which minimizes the quantity
\begin{equation*}
	\max_{x \in \mathcal N} \operatorname{dist}_\mathcal N(x,\phi_K(x)),
\end{equation*}
where $\operatorname{dist}_\mathcal N(\cdot,\cdot)$ denotes the graph distance, cf.~\cref{eq:dist}, and $\phi\colon \mathcal N \to C$ denoted the mapping to the nearest center with respect to $\operatorname{dist}_\mathcal N(\cdot,\cdot)$. The Gonzalez algorithm is a greedy-type algorithm that provides a 2-approximation to this problem with $\mathcal O(m \cdot \# \mathcal N)$ runtime. Interestingly, the algorithm is ``as good as possible'' in the sense that the existence of a $\delta$-approximation for $0<\delta<2$ would imply that $\mathrm{NP}$-hard problems can be solved in polynomial time, cf.~\cite{HoS85}.
After identifying the set of center nodes $C$, all remaining nodes are associated to the center nodes by the mapping $\phi$. An illustration of a spatial network partition generated by the Gonzalez algorithm is shown, e.g., in \cref{fig:partition}. We use an implementation of the Gonzalez algorithm from the toolbox GBFPUM, cf.~\cite{de2022gbfpum}.

\subsubsection*{Construction of a partition of unity}
Again, there are several ways to construct the partition of unity, and we will present only one possibility. In the following, we do not explicitly compute the mesh size $H$ defined in \cref{eq:H}, but use $H= N^{-1/d}$ instead, which is (up to a moderate factor) justified by \cref{ass:network} and the fact that $\Omega$ is scaled to unit size. We construct for any $T \in \TH$ the supports $U_T$ of the partition of unity function $\Lambda_T$ by extending the element~$T$ by all nodes with a graph distance to $T$ of less than or equal to $H/2$. For each $T \in \TH$, we then define~$\tilde \Lambda_T$ as $\tilde \Lambda_T(x)\coloneqq \operatorname{dist}_\mathcal N(x,\partial U_T)$ with $\partial U_T \coloneqq (\partial \Omega \cap U_T) \cup  \operatorname{bnd}_{\mathcal N}(U_T)$. If the element~$T$ includes boundary nodes on $\partial \Omega$, these nodes are excluded from $\partial U_T$. Finally, for all $T \in \TH$ the partition of unit function $\Lambda_T$ is obtained by
\begin{equation}
	\label{eq:pu}
	\Lambda_T \coloneqq \frac{\tilde \Lambda_T}{\sum_{K \in \mathcal T_H} \tilde \Lambda_K}.
\end{equation}
Since the overlap of the partition of unity functions is chosen sufficiently large, the denominator of \cref{eq:pu} is uniformly positive.

\begin{figure}
	\centering
	\includegraphics[height=0.325\linewidth]{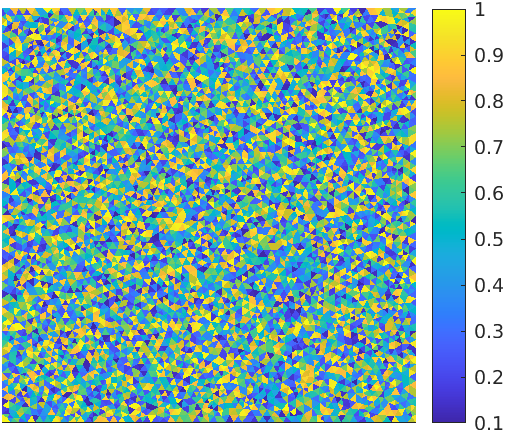}\hspace{1cm}
	\includegraphics[height=0.325\linewidth]{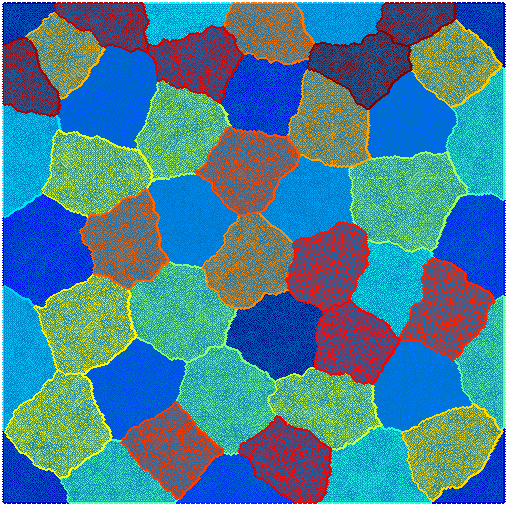}
	\caption{Piecewise constant coefficient $A$ (left) and partition of the finite element mesh with 50 elements (right).}
	\label{fig:partition_fem}
\end{figure}

\begin{samepage}
\subsection{Numerical examples}

\subsubsection*{Finite element discretization}

In this first numerical example, we apply our method to obtain a reliable coarse scale approximation of a fine finite element discretization. Let us consider the domain $\Omega = (0,1)^2$ and the weak formulation of the diffusion-type problem \end{samepage}

\begin{equation}
	\label{eq:diffprob}
	-\operatorname{div}(a \nabla u) = g_1
\end{equation}
subject to homogeneous Dirichlet boundary conditions. As source term we use $g_1(x,y) = \sin(x)\sin(y)$ and we choose a coefficient $a$ which is piecewise constant on a Delaunay triangulation $\mathcal T_a$ of $\Omega$ with around 8,500 elements. The element values of $a$ are chosen as realizations of independent random variables uniformly distributed in the interval $[0.1,1]$; see \cref{fig:partition_fem} (left) for an illustration of the coefficient. 
In order for the finite element method to resolve the multiscale features of the coefficient, we introduce the triangulation~$\mathcal T_h$ by twice global red refinement of $\mathcal T_a$.
We then apply our method for the coarse scale approximation of the finite element system $Ku = f$, where $K$ is the $\mathcal P^1$-finite element stiffness matrix and $f$ is the product of the $\mathcal P^1$-finite element mass matrix with the vector of nodal values of $g_1$.  To fit into our framework, we choose $\alpha = 0$, cf.~\cref{sec:choiceK}. The physical dimension $n = 2$ equals the isoperimetric dimension $d$ here.

\begin{figure}
	\centering
	\includegraphics[height=0.415\linewidth]{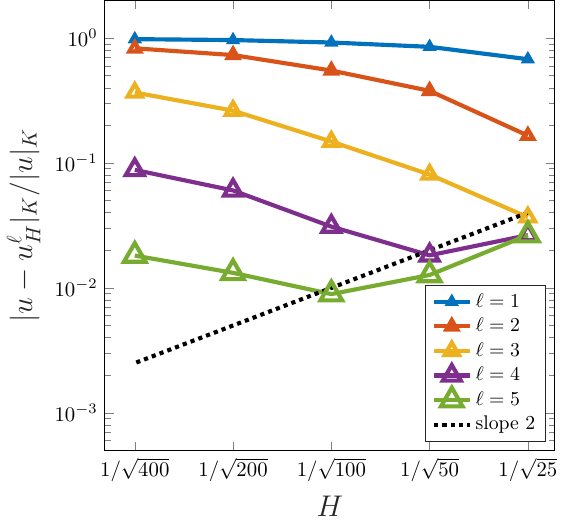}
	\hspace{3ex}
	\includegraphics[height=0.415\linewidth]{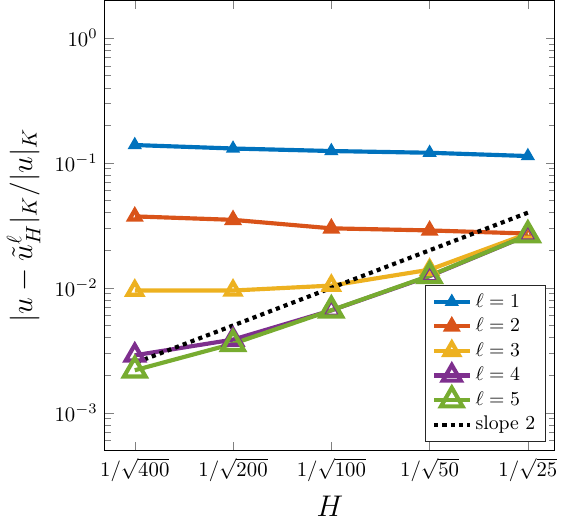}
	\caption{Relative errors for the naive (left) and stabilized (right) variants of the method applied to the finite element discretization for the source term $g_1$. The errors are plotted as functions of $H$ for several choices of~$\ell$.}
	\label{fig:conv_fem}
\end{figure}

We consider a hierarchy of partitions with $N\in \{25,50,100,200,400 \}$ elements, respectively; see \cref{fig:partition_fem} (right) for an illustration with $N = 50$.
In the following, we set $H = 1/\sqrt{N}$ which is, up to a moderate factor, equivalent to~\cref{eq:H}. This can be shown using  \cref{ass:network} and the fact that $\Omega$ is of unit size. In~\cref{fig:conv_fem}, we show the relative errors in the $K$-norm as a function of~$H$ for different values of the localization parameters~$\ell$, for the naive (left) and the stabilized (right) versions of the method. For the naive version, second-order convergence can be observed only for large localization parameters~$\ell$. For fixed $\ell$ and decreasing~$H$ one observes a deterioration of the error. The second-order convergence and the deterioration of the error for fixed $\ell$ are in agreement with \cref{thm:convsimple}. The stabilized version overcomes the latter effect. In particular, the second-order convergence can be seen much more clearly. For fixed localization parameters $\ell$, the error does not increase, but only stagnates. This is consistent with \cref{thm:convsimple}.

\subsubsection*{Planar spatial network model}

\begin{figure}
	\centering
	\includegraphics[height=0.325\linewidth]{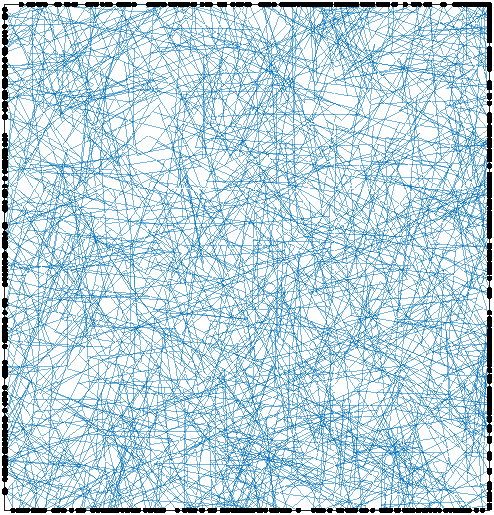}\hspace{1cm}
	\includegraphics[height=0.325\linewidth]{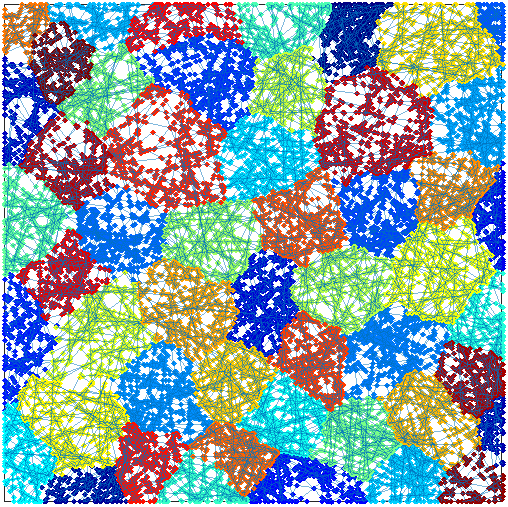}
	\caption{Spatial network with boundary nodes marked black (left) and partition of the spatial network with 50 elements (right).}
	\label{fig:partition}
\end{figure}

In the second numerical example, we consider a planar spatial network on $\Omega = (0,1)^2$. Its construction involves the following steps: First, we sample 4,000 lines of length 0.1, which are uniformly rotated and whose midpoints are uniformly distributed in $(-0.05,1.05)^2$. Next, we remove all line segments outside of $\Omega$. We then define the network nodes as the endpoints and intersections of the line segments. Two nodes are connected by an edge if they share a line segment. To ensure that the network is connected, we consider only the largest connected component of the network. We also remove all hanging nodes (nodes of degree one) that are not at the boundary $\partial \Omega$ along with their respective edges. All nodes on $\partial \Omega$ are Dirichlet nodes. The total number of nodes is about $144,000$. For an illustration of a similar spatial network, where we  use only 2,000 lines for illustration purposes, see~\cref{fig:partition} (left). Also for this spatial network the physical dimension $n = 2$ equals the isoperimetric dimension $d$. 
As operator~$K$, we consider the weighted graph Laplacian with parameter $\alpha = 1$, cf.~\cref{sec:choiceK}, and edge weights obtained as realizations of independent random variables which are uniformly distributed in the interval $[0.1,1]$. Due to the choice of edge weights, the problem exhibits multiscale features. In the following, we consider two problems of the form $Ku = f$ for two different right-hand sides. The first right-hand side is obtained by the multiplication of the diagonal mass matrix $M$ (cf.~\cref{eq:mass}) with the vector of nodal values of the source term $g_1$ from the previous numerical example. The second right-hand side is obtained similarly but with the function $g_2(x,y) \equiv 1$.

We consider the problem for the right-hand side involving $g_1$. Similarly as in the previous numerical example, we construct a hierarchy of partitions with $N\in \{25,50,100,200,400 \}$ elements and set $H = 1/\sqrt{N}$; see \cref{fig:partition} (right) for the visualization of a partition with $N = 50$.

\begin{figure}
	\centering
	\includegraphics[height=0.415\linewidth]{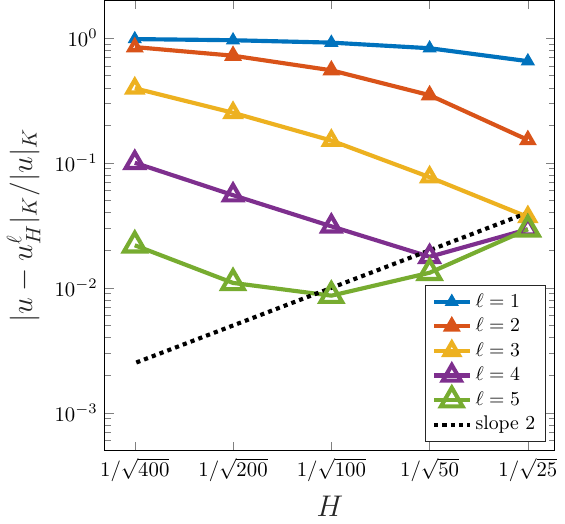}
	\hspace{3ex}
	\includegraphics[height=0.415\linewidth]{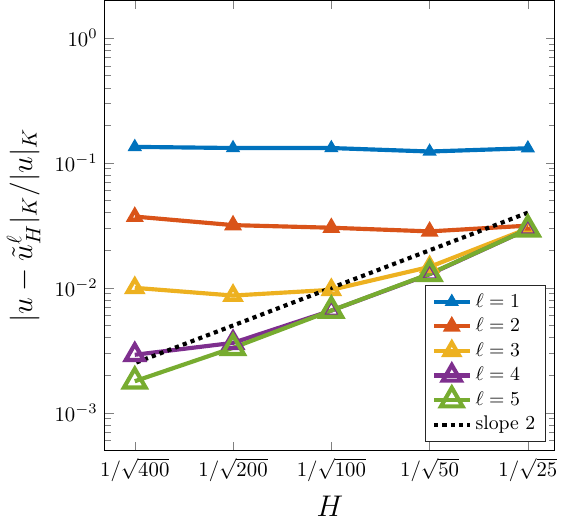}
	\caption{Relative errors for the naive (left) and stabilized (right) variants of the method applied to the planar spatial network for the source term $g_1$. The errors are plotted as functions of $H$ for several choices of~$\ell$.}
	\label{fig:conv}
\end{figure}

\cref{fig:conv} shows the relative errors in the $K$-norm as a function of~$H$ for different values of the localization parameters~$\ell$, for the naive (left) and the stabilized (right) versions of the method. Note that the qualitative behavior is almost the same as in the previous numerical example, which supports \cref{thm:convsimple,thm:convstab}.

\begin{figure}
	\centering\hspace{-1.75ex}
	\includegraphics[height=0.415\linewidth]{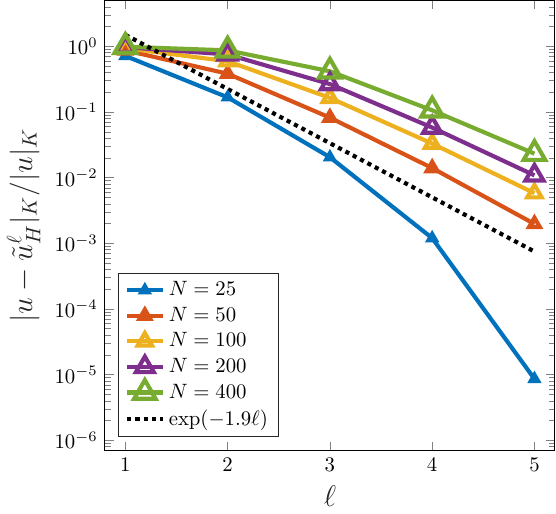}
	\hspace{4.5ex}
	\includegraphics[height=0.415\linewidth]{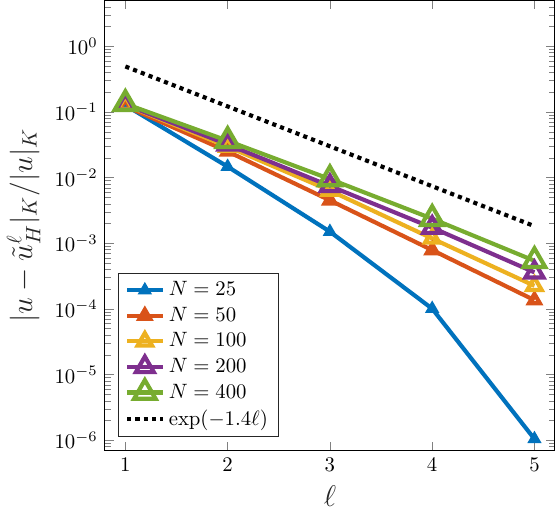}
	\caption{Relative errors for the naive (left) and stabilized (right) variants of the method applied to the planar spatial network for the source term $g_2$. The errors are plotted as functions of $\ell$ for several values of $H$.}
	\label{fig:loc}
\end{figure}

Second, we consider the problem involving the source term $g_2$. For this source, the optimal order term in the error estimates of \cref{thm:convsimple,thm:convstab} vanishes, leaving only the localization error, cf.~\cref{rem:exactnessideal}. We plot the corresponding error as a function of $\ell$ for several choices of $H$. In \cref{fig:loc} one observes an exponential decay of the localization error for the naive version (left) and the stabilized version (right), which is in line with~\cref{thm:convsimple,thm:convstab}. Note that the very fast decay for $N = 25$ is due to the fact that many patches coincide with the entire domain. For the naive version of the method, one observes that the error curves for different $\ell$ are shifted by a constant factor with respect to the previous one due to the negatively scaling term mentioned above. This is not observed for the stabilized version.

\subsubsection*{Cardboard-like spatial network model}

\begin{figure}
	\centering
	\includegraphics[height=0.3\linewidth]{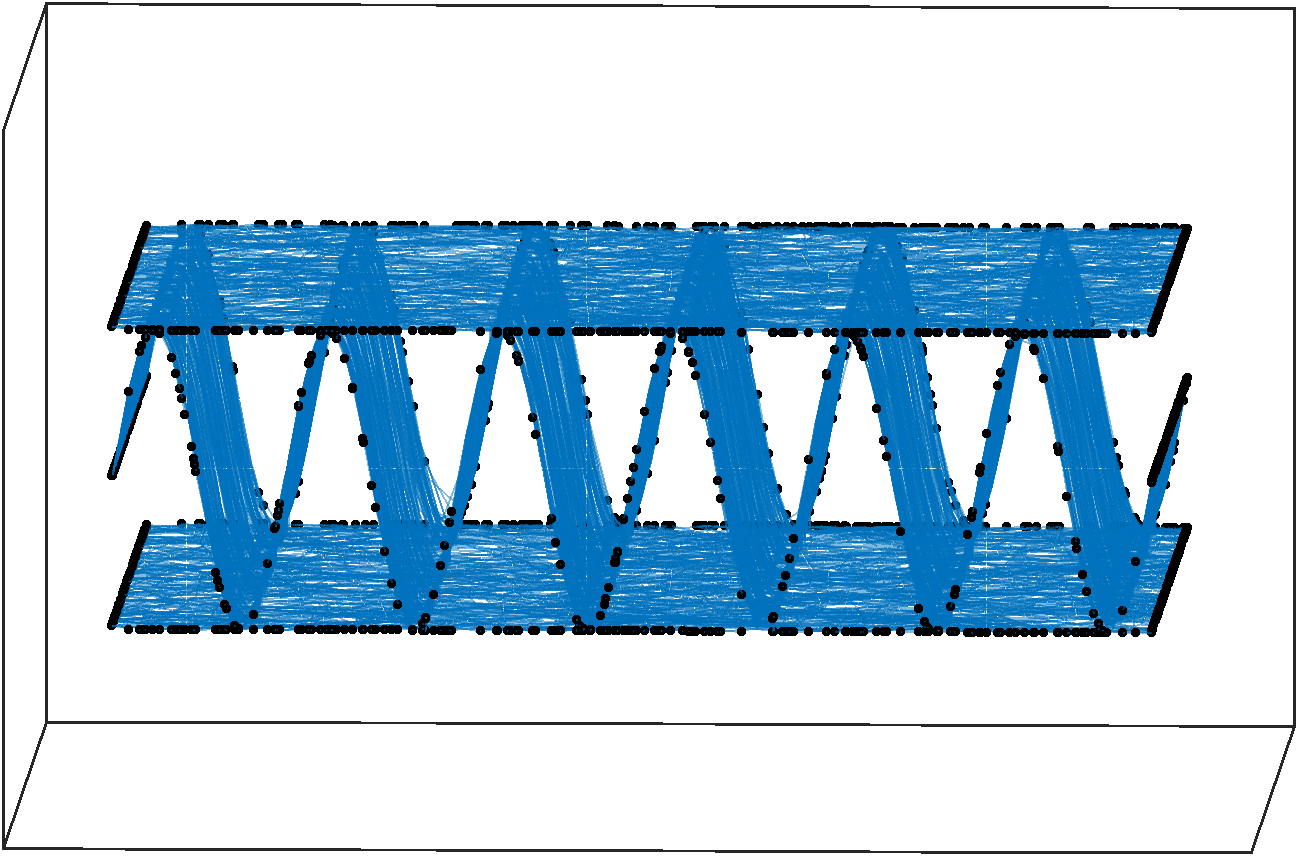}\hspace{2ex}
	\includegraphics[height=0.3\linewidth]{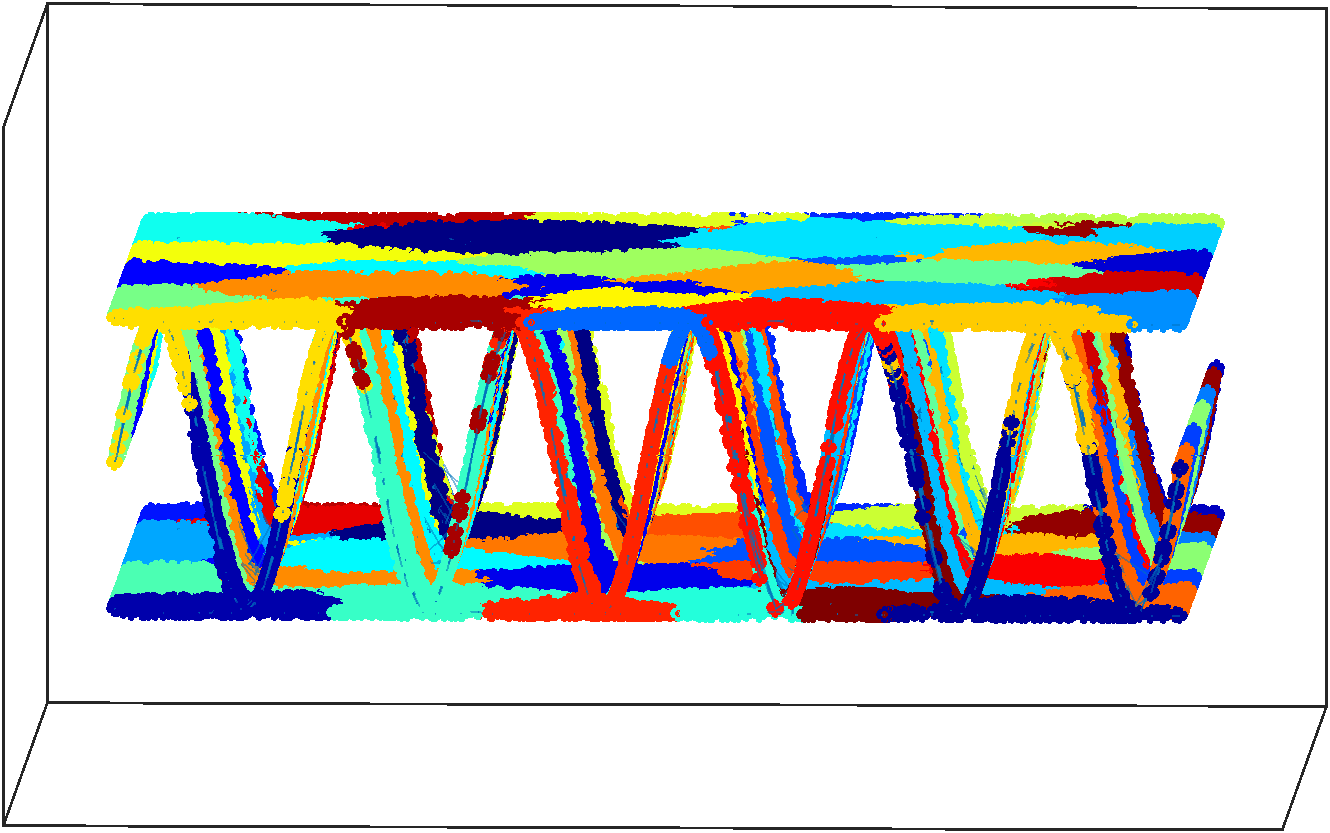}
	\caption{Cardboard-like spatial network with boundary nodes marked black (left) and illustration of a partition of this spatial network with 50 elements (right).}
	\label{fig:networkcarboard}
\end{figure}

Finally, we consider a spatial network that models a material that resembles corrugated cardboard. Such materials are of particular interest for the paper industry. We use the three-dimensional domain $\Omega = (0,1)^2\times (-\tfrac18-\delta,\tfrac18+\delta)$ with $\delta = 10^{-6}$ and construct the network in three steps: First, we create the top and bottom layers of the network using two planar networks with the same parameters as in the previous numerical example. For the top layer we set the $z$-coordinate of all nodes to $\tfrac18 + \delta$ and for the bottom layer to $-\tfrac18-\delta$. Second, we construct the oscillating layer between the top and bottom layers by adding the $z$ coordinate $z(x,y) = \tfrac18 \cos(12\pi x)$ to a planar spatial network again with the same parameters as in the previous numerical example. The last step is to sufficiently connect the three separate networks. To do this, we connect all nodes in the top and bottom layers to their nearest node in the oscillating layer provided that the distance is less than or equal to $10^{-4}$. Note that by altering $\delta$ one can control the minimal length of the edges connecting the layers. As boundary nodes we use the points whose $(x,y)$-coordinates lie on the boundary of the unit square. Note that this network has an isoperimetric dimension of $d=2$, cf.~\cref{ass:isoperimetric}, while being embedded into a domain of dimension $n=3$. The constructed cardboard-like spatial network has $433,000$ nodes.
For an illustration of a similar spatial network with fewer lines (for illustration purposes), see \cref{fig:networkcarboard} (left). As in the previous numerical example, we consider a weighted graph Laplacian with parameter $\alpha = 1$ and random edge weights in $[0.1,1]$ as the operator $K$ and use two right-hand sides for the source terms $(x,y,z)\mapsto g_i(x,y)$, $i \in \{1,2\}$. Note that $z$-dependent source term can also be handled, provided that the considered partition has sufficient resolution in the $z$-direction.

\begin{figure}
	\centering
	\includegraphics[height=.415\linewidth]{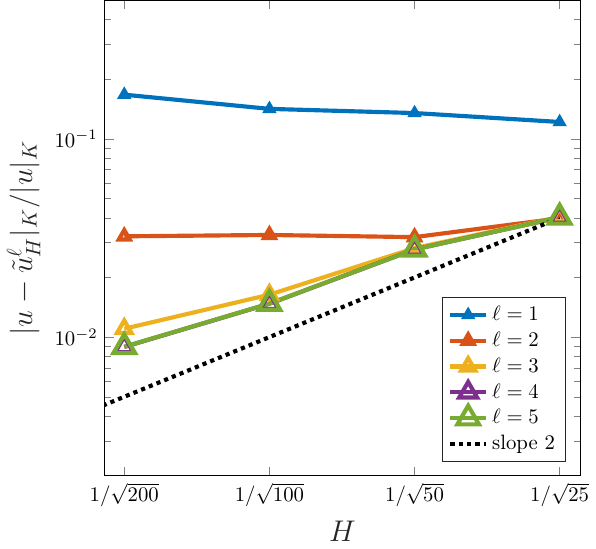}
	\hspace{4ex}
		\includegraphics[height=.415\linewidth]{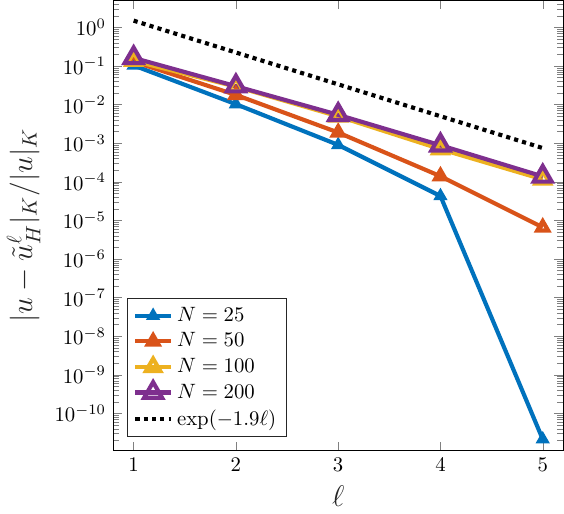}
	\caption{
Relative errors for the stable variant of the method applied to the cardboard-like  spatial network for the source terms $g_1$ (left) and~$g_2$~(right). The errors are plotted as functions of $H$ for several $\ell$ (left) and as functions of~$\ell$ for several $H$ (right).}
	\label{fig:stabilizedcardboard}
\end{figure}

In the following, we consider only the stabilized variant of the proposed method. To investigate its convergence, we use the right-hand side corresponding to the source term $g_1$ and use a hierarchy of partitions with $N \in \{25,50,100,200\}$ elements; see \cref{fig:networkcarboard} (right) for the visualization of one partition.
Since the isoperimetric dimension of the cardboard-like spatial network is two, we again set $H = 1/\sqrt{N}$. The plot of the relative $K$-norm errors is shown in \cref{fig:stabilizedcardboard}~(left). One observes a similar convergence as in the previous numerical experiments for the stabilized variant. For investigating the localization properties, we use the right-hand side corresponding to the source term $g_2$. In~\cref{fig:stabilizedcardboard} (right), the relative $K$-norm errors are plotted as functions of $\ell$ for several values of $H$. As before, one observes an exponential convergence of the error. Note that for $\ell = 5$ and $N = 25$ all patches are global which explains the exactness (modulo machine precision) of the method. We emphasize that also these numerical example clearly supports the theoretical results from \cref{thm:convstab}.

\section{Conclusion}

We have presented a multiscale method for spatial network models. We have shown that a system of equations which is spectrally equivalent to the graph Laplacian and has a similar sparsity pattern, can be reduced to a coarse-scale system of equations that suitably takes into account fine structures of the network. The methodology is built on the idea of the Localized Orthogonal Decomposition technique but works in an algebraic setting based on appropriate graph partitions. Our setting includes -- but is not limited to -- systems of equations resulting from classical finite element or finite difference discretizations of, e.g., elliptic partial differential equations. We have presented a full error analysis along with numerical examples. In particular, we have included an example with a complex material structure related to corrugated cardboard, which we could not consider before.

\bibliographystyle{alpha}
\bibliography{bib}

\newcommand{\etalchar}[1]{$^{#1}$}
\def\cprime{$'$}
\begin{thebibliography}{KMM{\etalchar{+}}20}

\bibitem[AHP21]{Peterseim2021}
R.~Altmann, P.~Henning, and D.~Peterseim.
\newblock Numerical homogenization beyond scale separation.
\newblock {\em Acta Numer.}, 30:1--86, 2021.

\bibitem[CEPT12]{CEPT12}
J.~Chu, B.~Engquist, M.~Prodanovi\'{c}, and R.~Tsai.
\newblock A multiscale method coupling network and continuum models in porous
  media i: steady-state single phase flow.
\newblock {\em Multiscale Model. Simul.}, 10:515--549, 2012.

\bibitem[Chu05]{C05}
F.~R.~K. Chung.
\newblock Laplacians and the cheeger inequality for directed graphs.
\newblock {\em Ann. Comb.}, 9:1--19, 2005.

\bibitem[CY95]{CY95}
F.~R.~K. Chung and S.-T. Yau.
\newblock Eigenvalues of graphs and {S}obolev inequalities.
\newblock {\em Comb. Probab. Comput.}, 4(1):11--25, 1995.

\bibitem[DHM23]{DHM23}
Z.~Dong, M.~Hauck, and R.~Maier.
\newblock An improved high-order method for elliptic multiscale problems.
\newblock {\em SIAM J. Numer. Anal.}, 61(4):1918--1937, 2023.

\bibitem[DRCE22]{de2022gbfpum}
A.~De~Rossi, R.~Cavoretto, and W.~Erb.
\newblock {GBFPUM-A MATLAB} package for partition of unity based signal
  interpolation and approximation on graphs.
\newblock {\em Dolomites Res. Notes Approx.}, 15(DRNA Volume 15.2):25--34,
  2022.

\bibitem[EGH{\etalchar{+}}24]{EGHKM24}
F.~Edelvik, M.~G\"{o}rtz, F.~Hellman, G.~Kettil, and A.~M\aa~lqvist.
\newblock Numerical homogenization of spatial network models.
\newblock {\em Comput. Methods Appl. Mech. Engrg.}, 418:Paper No. 116593, 2024.

\bibitem[EIL{\etalchar{+}}09]{EILRW09}
R.~Ewing, O.~Iliev, R.~Lazarov, I.~Rybak, and J.~Willems.
\newblock A simplified method for upscaling composite materials with high
  contrast of the conductivity.
\newblock {\em SIAM J.\ Sci.\ Comput.}, 31:2568--2586, 2009.

\bibitem[FHKP24]{pumslod}
P.~Freese, M.~Hauck, T.~Keil, and D.~Peterseim.
\newblock A super-localized generalized finite element method.
\newblock {\em to appear in Numer. Math.}, 2024.

\bibitem[FKO{\etalchar{+}}22]{FKOWW22}
M.~Fritz, T.~K{\"o}ppl, J.~T. Oden, A.~Wagner, and B.~Wohlmuth.
\newblock 1d--0d--3d coupled model for simulating blood flow and transport
  processes in breast tissue.
\newblock {\em Int. J. Numer. Method Biomed. Eng.}, 2022.

\bibitem[GHM22]{GoHeMa22}
M.~G{\"o}rtz, F.~Hellman, and A.~M{\aa}lqvist.
\newblock Iterative solution of spatial network models by subspace
  decomposition.
\newblock ArXiv e-print 2207.07488, 2022.

\bibitem[Gon85]{Gon85}
T.~F. Gonzalez.
\newblock Clustering to minimize the maximum intercluster distance.
\newblock {\em Theor. Comput. Sci.}, 38:293--306, 1985.

\bibitem[HM22]{HaM22}
M.~Hauck and A.~M{\aa}lqvist.
\newblock Super-localization of spatial network models.
\newblock ArXiv e-print 2210.07860, 2022.

\bibitem[HMM24]{hellman2022wellposedness}
F.~Hellman, A.~Målqvist, and M.~Mosquera.
\newblock Well-posedness and finite element approximation of mixed dimensional
  partial differential equations.
\newblock {\em to appear in BIT}, 2024.

\bibitem[HP13]{HeP13}
P.~Henning and D.~Peterseim.
\newblock Oversampling for the multiscale finite element method.
\newblock {\em Multiscale Model. Simul.}, 11(4):1149--1175, 2013.

\bibitem[HP22]{HaPe21}
M.~Hauck and D.~Peterseim.
\newblock Multi-resolution localized orthogonal decomposition for {H}elmholtz
  problems.
\newblock {\em Multiscale Model. Simul.}, 20(2):657--684, 2022.

\bibitem[HP23]{HaPe21b}
M.~Hauck and D.~Peterseim.
\newblock Super-localization of elliptic multiscale problems.
\newblock {\em Math. Comp.}, 92(341):981--1003, 2023.

\bibitem[HS85]{HoS85}
D.~S. Hochbaum and D.~B. Shmoys.
\newblock A best possible heuristic for the k-center problem.
\newblock {\em Math. Oper. Res.}, 10(2):180–184, 1985.

\bibitem[HWZ21]{HWZ21}
X.~Hu, K.~Wu, and L.~T. Zikatanov.
\newblock A posteriori error estimates for multilevel methods for graph
  {L}aplacians.
\newblock {\em SIAM J. Sci.\ Comput.}, 43(5):S727--S742, 2021.

\bibitem[ILW10]{ILW10}
O.~Iliev, R.~Lazarov, and J.~Willems.
\newblock Fast numerical upscaling of heat equation for fibrous materials.
\newblock {\em Comput.~Visual.~Sci.}, 13:275--285, 2010.

\bibitem[KMM{\etalchar{+}}20]{KMM20}
G.~Kettil, A.~M{\aa}lqvist, A.~Mark, M.~Fredlund, K.~Wester, and F.~Edelvik.
\newblock Numerical upscaling of discrete network models.
\newblock {\em BIT}, 60(1):67--92, 2020.

\bibitem[LB12]{LB12}
O.~E. Livne and A.~Brandt.
\newblock Lean algebraic multigrid (lamg): fast graph {L}aplacian linear
  solver.
\newblock {\em SIAM J. Sci.\ Comput.}, 34(4):B499--B522, 2012.

\bibitem[Mai21]{Mai20ppt}
R.~Maier.
\newblock A high-order approach to elliptic multiscale problems with general
  unstructured coefficients.
\newblock {\em SIAM J. Numer. Anal.}, 59(2):1067--1089, 2021.

\bibitem[MP14]{MaP14}
A.~M{\aa}lqvist and D.~Peterseim.
\newblock Localization of elliptic multiscale problems.
\newblock {\em Math. Comp.}, 83(290):2583--2603, 2014.

\bibitem[MP20]{MalP20}
A.~M{\aa}lqvist and D.~Peterseim.
\newblock {\em Numerical Homogenization by Localized Orthogonal Decomposition}.
\newblock Society for Industrial and Applied Mathematics (SIAM), Philadelphia,
  2020.

\bibitem[OS19]{OwhS19}
H.~Owhadi and C.~Scovel.
\newblock {\em Operator-Adapted Wavelets, Fast Solvers, and Numerical
  Homogenization: From a Game Theoretic Approach to Numerical Approximation and
  Algorithm Design}.
\newblock Cambridge Monographs on Applied and Computational Mathematics.
  Cambridge University Press, 2019.

\bibitem[Owh17]{Owh17}
H.~Owhadi.
\newblock Multigrid with rough coefficients and multiresolution operator
  decomposition from hierarchical information games.
\newblock {\em SIAM Rev.}, 59(1):99--149, 2017.

\bibitem[XZ17]{XZ17}
J.~Xu and L.~Zikatanov.
\newblock Algebraic multigrid methods.
\newblock {\em Acta Numer.}, 26:591--721, 2017.

\end{thebibliography}
\end{document}